\documentclass{amsart}
\usepackage{amssymb}

\usepackage{graphicx}
\usepackage{amscd}

\newtheorem{theorem}{Theorem}
\theoremstyle{plain}

\newtheorem{corollary}{Corollary}

\newtheorem{definition}{Definition}
\newtheorem{example}{Example}
\newtheorem{lemma}{Lemma}
\newtheorem{problem}{Problem}

\newtheorem{remark}{Remark}

\numberwithin{equation}{section}

\input{tcilatex}

\begin{document}
\title[On sub $B$-convex Banach spaces]{On sub $B$-convex Banach spaces}
\author{Eugene V. Tokarev}
\address{B.E. Ukrecolan, 33-81, Iskrinskaya str., 61005, Kharkiv-5, Ukraine}
\email{tokarev@univer.kharkov.ua}
\dedicatory{Dedicated to the memory of S. Banach.}

\begin{abstract}
In the article is introduced a new class of Banach spaces that are called
\textit{sub }$B$\textit{-convex}. Namely, a Banach space $X$ is said to be $%
B $-convex if it may be represented as a direct sum $l_{1}\oplus W$, where $%
W $ is $B$-convex. It will be shown that any separable sub $B$-convex Banach
space $X$ may be almost isometricaly embedded \ in a separable Banach space $%
G_{X}$ of the same cotype as $X$, which has a series of properties. Namely,

(1). $G_{X}$\textit{\ is an approximate envelope, i.e. any separable Banach
space }$Y$\textit{, which is finitely representable in }$G_{X}$\textit{, for
every }$\varepsilon >0$\textit{\ may be }$(1+\varepsilon )$\textit{%
-isomorphicaly embedded into }$G_{X}$;

(2). $G_{X}$\textit{\ is almost isotropic (i.e. has an almost transitive
norm)};

(3). $G_{X}$\textit{\ is existentialy closed in a class }$\left(
G_{X}\right) ^{f}$\textit{\ of all spaces that are finitely equivalent to }$%
G_{X}$;

(4). \textit{The conjugate space }$\left( G_{X}\right) ^{\ast }$\textit{\ is
of cotype 2};

(5). \textit{Every operator }$u:\left( G_{X}\right) ^{\ast }\rightarrow
l_{2} $\textit{is absolutely summing};

(6). \textit{Every projection }$P:G_{X}\rightarrow G_{X}$\textit{\ of finite
rank }$rank(P)=dim(PG_{X})$\textit{\ has a norm that is estimated as }$%
\left\| P\right\| \geq C\left( rank(P)\right) ^{1/q(X)}$\textit{, where }$%
q(X)$\textit{\ is the infimum of such }$q^{\prime }$\textit{s that }$X$%
\textit{\ is of cotype }$q$.

If, in addition, $X$ is of cotype 2, then $G_{X}$ has the more impressive
properties:

(7). \textit{Every operator }$v:G_{X}\rightarrow l_{2}$\textit{\ is
absolutely summing};

(8). \textit{All operators }$w:G_{X}\rightarrow l_{1}$\textit{\ and }$%
s:\left( G_{X}\right) ^{\ast }\rightarrow l_{1}$\textit{\ are absolutely
summing};

(9). \textit{An injective and a projective tensor products of }$G_{X}$%
\textit{\ by itself are identical}:
\begin{equation*}
G_{X}\overset{\wedge }{\otimes }G_{X}=G_{X}\overset{\vee }{\otimes }G_{X}.
\end{equation*}
\end{abstract}

\maketitle

\section{Introduction}

Recall that a Banach space $X$ is $B$\textit{-convex }provided the space $%
l_{1}$ is not finitely representable in $X$ (exact definitions of all
notions will be given below).

In this article is introduced a new class of Banach spaces that are called
\textit{sub }$B$\textit{-convex}. Namely, a Banach space $X$ is said to be $%
B $-convex if it may be represented as a direct sum $l_{1}\oplus W$, where $%
W $ is $B$-convex.

It will be shown that any separable sub $B$-convex Banach space $X$ may be
almost isometricaly embedded \ in a separable Banach space $G_{X}$ of the
same cotype as $X$, which has a series of properties, some of which are
rather peculiar. Namely,

\begin{enumerate}
\item  $G_{X}$\textit{\ is an approximate envelope, i.e. any separable
Banach space }$Y$\textit{, which is finitely representable in }$G_{X}$%
\textit{, for every }$\varepsilon >0$\textit{\ may be }$(1+\varepsilon )$%
\textit{-isomorphicaly embedded into }$G_{X}$\textit{.}

\item  $G_{X}$\textit{\ is almost isotropic (i.e. has an almost transitive
norm).}

\item  $G_{X}$\textit{\ is existentialy closed in a class }$\left(
G_{X}\right) ^{f}$\textit{\ of all spaces that are finitely equivalent to }$%
G_{X}.$

\item  \textit{The conjugate space }$\left( G_{X}\right) ^{\ast }$\textit{\
is of cotype 2.}

\item  \textit{Every operator }$u:$\textit{\ }$\left( G_{X}\right) ^{\ast
}\rightarrow l_{2}$\textit{\ is absolutely summing.}

\item  \textit{Every projection }$P:G_{X}\rightarrow G_{X}$\textit{\ of
finite rank }$rank(P)=dim(PG_{X})$\textit{\ has a norm that is estimated as }%
$\left\| P\right\| \geq C\left( rank(P)\right) ^{1/q(X)}$\textit{, where }$%
q(X)$\textit{\ is the infimum if such }$q^{\prime }$\textit{s that }$X$%
\textit{\ is of cotype }$q$\textit{.}

If, in addition, $X$ is of cotype 2, then $G_{X}$ has the more impressive
properties:

\item  \textit{Every operator }$v:G_{X}\rightarrow l_{2}$\textit{\ is
absolutely summing.}

\item  \textit{All operators }$w:G_{X}\rightarrow l_{1}$\textit{\ and }$%
s:\left( G_{X}\right) ^{\ast }\rightarrow l_{1}$\textit{\ are absolutely
summing.}

\item  \textit{An injective and a projective tensor products of }$G_{X}$%
\textit{\ by itself are identical:}
\begin{equation*}
G_{X}\overset{\wedge }{\otimes }G_{X}=G_{X}\overset{\vee }{\otimes }G_{X}.
\end{equation*}
\end{enumerate}

The idea of construction is based on a notion of finite equivalence: Banach
spaces $X$ and $Y$ are\textit{\ finitely equivalent} ($X\sim _{f}Y$)
provided $X$ is finitely representable in $Y$ and $Y$ is finitely equivalent
to $X$ too. Thus, any Banach space $X$ generates a class $X^{f}=\{Y:X\sim
_{f}Y\}$ and a set $\frak{M}\left( X^{f}\right) $ of all finite dimensional
Banach spaces that are finitely representable in $X$. Obviously, $\frak{M}%
\left( X^{f}\right) \subseteq \frak{M}\left( Y^{f}\right) $ if and only if $%
X $ is finitely representable in $Y$.

The main step in construction is to find an enlargement $\frak{M}\left(
Y^{f}\right) $ of $\frak{M}\left( X^{f}\right) $, where $X$ is a given sub $%
B $-convex Banach space, in such a way that $Y$ will be of the same cotype
as $X$, and $\frak{M}\left( Y^{f}\right) $ has the so called \textit{%
amalgamation property}. This means that for any fifth $\left\langle
A,B_{1},B_{2},i_{1},i_{2}\right\rangle $, where $A$, $B_{1}$, $B_{2}\in
\frak{M}(X^{f})$; $i_{1}:A\rightarrow B_{1}$ and$\ i_{2}:A\rightarrow B_{2}$
are isometric embedding, there exists a triple $\left\langle
j_{1},j_{2},F\right\rangle $, where $F\in \frak{M}(X^{f})$; $\
j_{1}:B_{1}\rightarrow F$ and $j_{2}:B_{2}\rightarrow F$ are isometric
embedding such that $j_{1}\circ i_{1}=j_{2}\circ i_{2}$.

This property give a chance to construct a separable space $G_{X}$, which
has properties 1, 2 and 3, listed above.

The proof of other properties of $G_{X}$ is based on J. Pisier's results
[1], where the first example of a Banach space that has properties 4, 5, 7,
8, 9 and the property 6 in a case when $G_{X}$ is of cotype 2 was
constructed.

\section{Definitions and notations}

\begin{definition}
Let $X$, $Y$ are Banach spaces. $X$ is \textit{finitely representable} in $Y$
(in symbols: $X<_{f}Y$) if for each $\varepsilon >0$ and for every finite
dimensional subspace $A$ of $X$ there exists a subspace $B$ of $Y$ and an
isomorphism $u:A\rightarrow B$ such that $\left\| u\right\| \left\|
u^{-1}\right\| \leq 1+\varepsilon $.

Spaces $X$ and $\ Y$ are said to be finitely equivalent if $X<_{f}Y$ and $%
Y<_{f}X$.

Any Banach space $X$ generates classes
\begin{equation*}
X^{f}=\{Y\in \mathcal{B}:X\sim _{f}Y\}\text{ \ and \ }X^{<f}=\{Y\in \mathcal{%
B}:Y<_{f}X\}
\end{equation*}
\end{definition}

For any two Banach spaces $X$, $Y$ their \textit{Banach-Mazur distance }is
given by
\begin{equation*}
d(X,Y)=\inf \{\left\| u\right\| \left\| u^{-1}\right\| :u:X\rightarrow Y\},
\end{equation*}
where $u$ runs all isomorphisms between $X$ and $Y$ and is assumed, as
usual, that $\inf \varnothing =\infty $.

It is well known that $\log d(X,Y)$ defines a metric on each class of
isomorphic Banach spaces. A set $\frak{M}_{n}$ of all $n$-dimensional Banach
spaces, equipped with this metric, is a compact metric space that is called
\textit{the Minkowski compact} $\frak{M}_{n}$.

A disjoint union $\cup \{\frak{M}_{n}:n<\infty \}=\frak{M}$ is a separable
metric space, which is called the \textit{Minkowski space}.

Consider a Banach space $X$. Let $H\left( X\right) $ be a set of all its
different finite dimensional subspaces (isometric finite dimensional
subspaces of $X$ in $H\left( X\right) $ are identified). Thus, $H\left(
X\right) $ may be regarded as a subset of $\frak{M}$, equipped with the
restriction of the metric topology of $\frak{M}$.

Of course, $H\left( X\right) $ need not to be a closed subset of $\frak{M}$.
Its closure in $\frak{M}$ will be denoted $\overline{H\left( X\right) }$.
From definitions follows that $X<_{f}Y$ if and only if $\overline{H\left(
X\right) }\subseteq \overline{H\left( Y\right) }$. Spaces $X$ and $Y$ are
\textit{finitely equivalent }(in symbols: $X\sim _{f}Y$) if simultaneously $%
X<_{f}Y$ and $Y<_{f}X$. Therefore, $X\sim _{f}Y$ if and only if $\overline{%
H\left( X\right) }=\overline{H\left( Y\right) }$.

There exists a one to one correspondence between classes of finite
equivalence $X^{f}=\{Y\in \mathcal{B}:X\sim _{f}Y\}$ and closed subsets of $%
\frak{M}$ of kind $\overline{H\left( X\right) }$.

Indeed, all spaces $Y$ from $X^{f}$ have the same set $\overline{H\left(
X\right) }$. This set, uniquely determined by $X$ (or, equivalently, by $%
X^{f}$), will be denoted by $\frak{M}(X^{f})$ and will be referred as to
\textit{the Minkowski's base of the class} $X^{f}$.

Using this correspondence, it may be defined a set $f\left( \mathcal{B}%
\right) $ of all different classes of finite equivalence, assuming (to
exclude contradictions with the set theory) that members of $f\left(
\mathcal{B}\right) $ are sets $\frak{M}(X^{f})$. For simplicity it may be
says that members of $f\left( \mathcal{B}\right) $ are classes $X^{f}$
itself.

Clearly, $f\left( \mathcal{B}\right) $ is partially ordered by the relation $%
\frak{M}(X^{f})\subseteq \frak{M(}Y^{f})$, which may be replaced by the
relation $X^{f}<_{f}Y^{f}$ of the same meaning. The minimal (with respect to
this order) element of $\ f\left( \mathcal{B}\right) $ is the class $\left(
l_{2}\right) ^{f}$ (the Dvoretzki theorem); the maximal one - the class $%
\left( l_{\infty }\right) ^{f}$ (an easy consequence of the Hahn-Banach
theorem). Other $l_{p}$'s are used in the classifications of Banach spaces,
which was proposed by L. Schwartz [2].

For a Banach space $X$ its $l_{p}$-\textit{spectrum }$S(X)$ is given by
\begin{equation*}
S(X)=\{p\in\lbrack0,\infty]:l_{p}<_{f}X\}.
\end{equation*}

Certainly, if $X\sim_{f}Y$ then $S(X)=S(Y)$. Thus, the $l_{p}$-spectrum $%
S(X) $ may be regarded as a property of the whole class $X^{f}$. So,
notations like $S(X^{f})$ are of obvious meaning.

Let $X$ be a Banach space. It is called:

\begin{itemize}
\item  $c$-\textit{convex,} if $\infty \notin S(X)$;

\item  $B$-\textit{convex,} if $1\notin S\left( X\right) $;

\item  \textit{Finite universal,} if $\infty \in S(X)$.

\item  \textit{Superreflexive,} if every space of the class $X^{f}$ is
reflexive.
\end{itemize}

Equivalently, $X$ is superreflexive if any $Y<_{f}X$ is reflexive. Clearly,
any superreflexive Banach space is $B$-convex.

As was shown in [3], the $l_{p}$-spectrum is closely connected with notions
of type and cotype. Recall the definition.

Let $1\leq p\leq2\leq q\leq\infty$. A Banach space $X$ is said to be of
\textit{type} $p$ (resp., of \textit{cotype }$q$) if for every finite
sequence $\{x_{n}:n<N\}$ of elements of $X$%
\begin{equation*}
\int_{0}^{1}\left\| \sum_{n=0}^{N-1}r_{n}\left( t\right) x_{n}\right\|
dt\leq t_{p}\left( X\right) \left( \sum_{n=0}^{N-1}\left\| x_{n}\right\|
^{p}\right) ^{1/p}
\end{equation*}
(respectively,
\begin{equation*}
\left( \sum_{n=0}^{N-1}\left\| x_{n}\right\| ^{q}\right) ^{1/q}\leq
c_{q}\left( X\right) \int_{0}^{1}\left\| \sum_{n=0}^{N-1}r_{n}\left(
t\right) x_{n}\right\| dt),
\end{equation*}
where $\left\{ r_{n}\left( t\right) :n<\infty\right\} $ are Rademacher
functions.

When $\ q=\infty$, the sum $\left( \sum_{n=0}^{N-1}\left\| x_{n}\right\|
^{q}\right) ^{1/q}$ must be replaced with $\underset{n<N}{\sup}\left\|
x_{n}\right\| $. Constants $t_{p}\left( X\right) $ and $c_{q}\left( X\right)
$\ in these inequalities depend only on $X$. Their least meanings, $%
T_{p}\left( X\right) $ and $C_{q}\left( X\right) $ respectively, are named
\textit{the type} $p${\small -}\textit{constant} $T_{p}\left( X\right) $,
resp., \textit{the cotype} $q$\textit{-constant }$C_{q}\left( X\right) $.

Any Banach space is of type $1$ and of cotype $\infty$.

If $X$ is of type $p$ and of cotype $q$ with the constants $T_{p}\left(
X\right) =T$, $C_{q}\left( X\right) =C$, than any $Y\in X^{f}$ is of same
type and cotype and its type-cotype constants are equal to that of $X$.
Thus, we may speak about type and cotype of the whole class $X^{f}$.

It is known (see [3]) that
\begin{align*}
\inf S(X)& =\sup \{p\in \lbrack 1,2]:T_{p}\left( X\right) <\infty \}; \\
\sup S(X)& =\inf \{q\in \lbrack 2,\infty ]:C_{q}\left( X\right) <\infty \}.
\end{align*}

In what follows a notion of inductive (or direct) limit will be used. Recall
a definition.

Let $\left\langle I,\ll\right\rangle $ be a \textit{partially ordered set}.
It said to be \textit{directed} (to the right hand) if for any $i,j\in I$
there exists $k\in I$ such that $i\ll k$ and $j\ll k$ .

Let $\left\{ X_{i}:i\in I\right\} $ be a set of Banach spaces that are
indexed by elements of an directed set $\left\langle I,\ll\right\rangle $.
Let $m_{i,j}:X_{i}\rightarrow X_{j}$ $(i\ll j)$ be isometric embedding.

A system $\left\{ X_{i},m_{i,j}:i,j\in I;i\ll j\right\} $ is said to be an
\textit{inductive (}or\textit{\ direct) isometric system} if
\begin{equation*}
m_{i,i}=Id_{X_{i}};\text{ \ \ }m_{i,k}=m_{j,k}\cdot m_{i,j}
\end{equation*}
for all $i\ll j\ll k$ ($Id_{Y}$ denotes the identical operator on $Y$).

Its \textit{inductive }(or \textit{direct}) \textit{limit}, $\underset{%
\rightarrow }{{\lim }}X_{i}$ is defined as follows. Let
\begin{equation*}
X=\cup \left\{ X_{i}\times \left\{ i\right\} :i\in I\right\}
\end{equation*}

Elements of $X$ are pairs $(x,i),$ where $x\in X_{i}$. Let $=_{eq}$ be a
relation of equivalence of elements of $X$, which is given by the following
rule:

\begin{center}
$(x,i)=_{eq}(y,j)$ if $m_{i,k}x=m_{j,k}y$ for some $k\in I$.
\end{center}

A class of all elements of $X$ that are equivalent to a given $(x,i)$ is
denoted as
\begin{equation*}
\lbrack x,i]=\{(y,j):(y,j)=_{eq}(x,i)\}.
\end{equation*}

A set of all equivalence classes $[x,i]$ is denoted $X_{\infty}$. Clearly, $%
X_{\infty}$ is a linear space. Let $\left\| [x,i]\right\| =\lim_{I}\left\|
m_{i,j}x\right\| _{X_{j}}$ be a semi-norm on $X_{\infty}$. Let
\begin{equation*}
Null(X)=\{[x,i]:\left\| [x,i]\right\| =0\}.
\end{equation*}

A direct limit of the inductive isometric system $\left\{
X_{i},m_{i,j}:i,j\in I;i\ll j\right\} $ is a quotient space
\begin{equation*}
\underset{\rightarrow}{\lim}X_{i}=X_{\infty}/Null(X).
\end{equation*}

Clearly, if all spaces $X_{i}$ ($i\in I$) belong to a given class $Z^{f}$,
then $\underset{\rightarrow }{\lim }X_{i}\in Z^{f}$; if all $X_{i}$'s are
finitely representable in $Z$ then $\underset{\rightarrow }{\lim }%
X_{i}<_{f}Z $ too.

\section{The amalgamation property and the space $l_{1}$}

It will be convenient to introduce some terminology.

A fifth $v=\left\langle A,B_{1},B_{2},i_{1},i_{2}\right\rangle $, where $A$,
$B_{1}$, $B_{2}\in \frak{M}(X^{f})$; $i_{1}:A\rightarrow B_{1}$ and$\
i_{2}:A\rightarrow B_{2}$ are isometric embedding will be called \textit{the
}$V$\textit{-formation over} $\frak{M}(X^{f})$. The space $A$ will be called
\textit{the root of the }$V$\textit{-formation }$v$\textit{. If there exists
}a triple $t=\left\langle j_{1},j_{2},F\right\rangle $, where $F\in \frak{M}%
(X^{f})$; $\ j_{1}:B_{1}\rightarrow F$ and $j_{2}:B_{2}\rightarrow F$ are
isometric embedding such that $j_{1}\circ i_{1}=j_{2}\circ i_{2}$, then the $%
V$-formation $v$ is said to be \textit{amalgamated in} $\frak{M}(X^{f})$,
and the triple $t$ is said to be its \textit{amalgam}.

Let $Amalg(\frak{M}(X^{f}))$ be a set of all spaces $A\in \frak{M}(X^{f})$
with the property:

\textit{Any }$V$\textit{-formation }$v$\textit{, which root is }$A$\textit{\
is amalgamated in }$\frak{M}(X^{f})$.

\begin{definition}
Let $X\in \mathcal{B}$ \ generates a class $X^{f}$\ with a Minkowski's base $%
\frak{M}(X^{f})$. It will be said that $\frak{M}(X^{f})$ (and the whole
class $X^{f}$) has the amalgamation property if
\begin{equation*}
\frak{M}(X^{f})=Amalg(\frak{M}(X^{f}))
\end{equation*}
\end{definition}

In a general setting of arbitrary algebraic structures the amalgamation
property was introduced by B. J\'{o}nsson [4]. To a Banach space case it was
transferred in the abstract of E.D. Positselski and author [5].

All results of this section were obtained in collaboration with E.D.
Positselski, whose tragic death terminated a study of this important
property.

To show that a set $\frak{M}(\left( l_{1}\right) ^{f})$ has the amalgamation
property it will be needed some additional constructions.

\begin{definition}
Let $X$ be a Banach space with a basis $\left( e_{n}\right) $. This basis is
said to be Besselian if there exists such a constant $C>0$ that for any $%
x=\sum x_{n}e_{n}\in X$%
\begin{equation*}
\left\| x\right\| \geq C\left( \sum \left( x_{n}\right) ^{2}\right) ^{1/2}.
\end{equation*}
\end{definition}

\begin{theorem}
Let $\left\langle X,\left( e_{n}\right) \right\rangle $ be a Banach space
with a Besselian bases $\left( e_{n}\right) $; $Y\hookrightarrow X$ be a
finite dimensional subspace of $X$ of dimension $dim(Y)=m$.

There exists a sequence of points $Y_{0}=\left( y_{n}\right) _{n=1}^{\infty
}\subset \mathbb{R}^{m}$ such that $X$ may be represented as a Banach space $%
\overline{X}$ of all rear-valued functions $f$, defined on $Y_{0}$ with a
norm
\begin{equation*}
\left\| f\right\| \overset{def}{=}\left\| f\right\| _{\overline{X}}=\left\|
\sum f\left( y_{i}\right) e_{i}\right\| _{X}<\infty ;
\end{equation*}
$Y$ may be represented as a subspace $\overline{Y}$ of $\overline{X}$ that
consists of all linear functions on $Y_{0}$ (equipped with the corresponding
restriction of the norm of $\overline{X}$). In other words, there exists an
isometry $T:X\rightarrow \overline{X}$ such that its restriction $T\mid
_{Y}:Y\rightarrow \overline{Y}$ is also an isometry.
\end{theorem}

\begin{proof}
Let $u=\sum u_{n}e_{n}$ and $v=\sum v_{n}e_{n}$ be elements of $X$. Since $%
\left( e_{n}\right) $ is Besselian, it may be defined a scalar product $%
\left\langle u,v\right\rangle =\sum u_{n}v_{n}$. Since $Y$is of finite
dimension, there exists a projection $P:X\rightarrow Y$ such that $%
\left\langle Px,y\right\rangle =\left\langle x,y\right\rangle $. Let $%
y_{i}=Pe_{i}$.

Consider a linear space $\overline{X}$ of all formal sums $f=\sum f_{i}y_{i}$
and equip it with a norm
\begin{equation*}
\left\| f\right\| _{\overline{X}}=\left\| \sum f_{i}y_{i}\right\| _{%
\overline{X}}=\left\| \sum f_{i}e_{i}\right\| _{X}.
\end{equation*}

Let a map $T:X\rightarrow \overline{X}$ be given by
\begin{equation*}
T(\sum u_{n}e_{n})=\sum u_{n}y_{n}\text{ \ for all }u=\sum u_{n}e_{n}\in X.
\end{equation*}

If $u\in Y$ then the $i$'s coordinate
\begin{equation*}
\left( Tu\right) _{i}=u_{i}=\left\langle u,e_{i}\right\rangle =\left\langle
u,Pe_{i}\right\rangle =\left\langle u,y_{i}\right\rangle .
\end{equation*}

Hence, $T$ maps $Y$ onto a linear space $\overline{Y}$ of all linear
functions, defined on $Y_{0}$.

So, $\overline{X}$, $\overline{Y}$, $Y_{0}$ and $T$ have desired properties.
\end{proof}

\begin{remark}
A set $Y_{0}=\left( y_{i}\right) \subset \mathbb{R}^{m}$ will be called an
incarnating set for a pair $Y\hookrightarrow X$. A pair $[\overline{X},%
\overline{Y}]$ will be called an incarnation pair for $Y\hookrightarrow X$.
\end{remark}

\begin{remark}
If a basis $\left( e_{n}\right) $ of $X$ is in addition a symmetric one, a
set $Y_{0}$ may be symmetrized:

it may be considered instead $Y_{0}$ a central symmetric set
\begin{equation*}
K=K\left( Y\hookrightarrow X\right) =Y_{0}\cup \left( -Y_{0}\right) =\{y\in
\mathbb{R}^{m}:y\in Y_{0}\text{ \ or \ }-y\in Y_{0}\}
\end{equation*}
as an incarnating set for a pair $Y\hookrightarrow X$, if as $\overline{X}\ $%
(resp., as $\overline{Y}$) will be considered all odd (resp., all odd
linear) functions on $K$ with the corresponding norm.
\end{remark}

Consider a space $\mathbb{R}^{m}$. Let $\left( e_{i}\right) _{i=1}^{m}$be a
basis of $\mathbb{R}^{m}$ that will be assumed to be orthogonal (with
respect to $l_{2}$-norm). Define on $\mathbb{R}^{m}$ an other norm, namely,
let
\begin{equation*}
\left| \left| \left| u\right| \right| \right| =\left| \left| \left|
\sum\nolimits_{i=1}^{m}u_{n}e_{n}\right| \right| \right|
=\sum\nolimits_{i=1}^{m}\left| \left\langle u,e_{i}\right\rangle \right|
=\sum\nolimits_{i=1}^{m}\left| u_{i}\right| .
\end{equation*}

\begin{theorem}
Let $m\in \mathbb{N}$, $K=\{\pm y_{i}:i\in \mathbb{N}\}\subset \mathbb{R}%
^{m} $ be a central symmetric set. $K$ is an incarnating set for a pair $%
Y\hookrightarrow l_{1}$ where $Y$ is $m$-dimensional subspace of $l_{1}$ if
and only if $K$ is complete in $\mathbb{R}^{m}$ and $\sum \{\left| \left|
\left| y_{i}\right| \right| \right| :y_{i}\in K\}<\infty $.
\end{theorem}

\begin{proof}
Certainly, the condition of completeness of $K$ in $\mathbb{R}^{m}$ is a
necessary one. The condition $\sum \{\left| \left| \left| y_{i}\right|
\right| \right| :y_{i}\in K\}<\infty $ is also necessary. Indeed,
\begin{eqnarray*}
\sum \left| \left| \left| y_{i}\right| \right| \right|
&=&\sum\nolimits_{y_{i}\in K}\sum\nolimits_{j=1}^{m}\left| \left\langle
y_{i},e_{j}\right\rangle \right| \\
&=&\sum\nolimits_{j=1}^{m}\sum\nolimits_{y_{i}\in K}\left| \left\langle
y_{i},e_{j}\right\rangle \right| =\sum\nolimits_{j=1}^{m}\left\|
e_{j}\right\| =m<\infty .
\end{eqnarray*}

The sufficiency follows from a fact that a set $L(K)$ of all linear
functionals on $K$ may be embedded in $l_{1}$. Indeed, let $x\in L(K)$; $%
x=\sum \{x_{i}y_{i}:y_{i}\in K\}$. Then
\begin{eqnarray*}
\left\| x\right\| _{l_{1}} &=&\sum\nolimits_{y_{i}\in K}\left| \left\langle
x,y_{i}\right\rangle \right| =\sum\nolimits_{y_{i}\in K}\left|
\sum\nolimits_{y_{j}\in K}x_{j}\left\langle y_{i},y_{j}\right\rangle \right|
\\
&\leq &\sum\nolimits_{y_{i}\in K}\max \left| x_{j}\right|
\sum\nolimits_{y_{j}\in K}\left| \left\langle y_{j},y_{i}\right\rangle
\right| =\max \left| x_{j}\right| \sum\nolimits_{y_{i}\in K}\left| \left|
\left| y_{i}\right| \right| \right| .
\end{eqnarray*}
\end{proof}

\begin{remark}
Obviously, if linearly congruent sets $K$ and $K_{^{1}}$ are incarnating
sets for pairs $Y\hookrightarrow l_{1}$ and $Z\hookrightarrow l_{1}$, then
the linear congruence $U:K\rightarrow K_{1}$ generates an isometric
automorphism $u:l_{1}\rightarrow l_{1}$ which restriction to $Y$ is an
isometry between $Y$ and $Z$.
\end{remark}

The following result shows how to reconstruct a pair $Y\hookrightarrow l_{1}$
by a given set $K=\{\pm y_{i}:i\in \mathbb{N}\}\subset \mathbb{R}^{m}$ that
is complete in $\mathbb{R}^{m}$ and $\sum \{\left| \left| \left|
y_{i}\right| \right| \right| :y_{i}\in K\}<\infty $.

Let
\begin{equation*}
\varsigma K=\cup \{\sum \{y_{i}\in K:\left\langle y_{i},u\right\rangle
>0\}:u\in \mathbb{R}^{m}\};\text{ \ }K^{\prime }=\overline{co\left(
\varsigma K\right) },
\end{equation*}
i.e. $K^{\prime }$ is a closure of a convex hull of $\varsigma K$.

\begin{theorem}
A central symmetric complete in $\mathbb{R}^{m}$ subset $K=\{\pm y_{i}:i\in
\mathbb{N}\}$ is an incarnating set for a pair $Y\hookrightarrow l_{1}$ if
and only if $K^{\prime }$ is congruent with the unit ball
\begin{equation*}
B(Y^{\ast })\overset{def}{=}\{y^{\prime }\in Y^{\ast }:\left\| y^{\prime
}\right\| \leq 1\}
\end{equation*}
of the conjugate space $Y^{\ast }$.
\end{theorem}

\begin{proof}
Necessity.
\begin{equation*}
\left\| x\right\| _{Y}=\frac{1}{2}\sum\nolimits_{y_{i}\in K}\left|
\left\langle x,y_{i}\right\rangle \right| =\max \{\left\langle
x,z\right\rangle :z\in K^{\prime }\}.
\end{equation*}
Hence, $K^{\prime }=B(Y^{\ast })$.

Sufficiency. $K$ is an incarnating set for a some pair $Z\hookrightarrow
l_{1}$. Hence, as was shown before, $K^{\prime }=B(Z^{\ast })$. However, by
conditions of the theorem, $K^{\prime }=B(Y^{\ast })$.

Thus, $Y^{\ast }=Z^{\ast }$ and, hence, $Y=Z$.
\end{proof}

\begin{corollary}
There exists the unambiguous correspondence between pairs $Z\hookrightarrow
l_{1}$ and certain central symmetric complete subsets of $\mathbb{R}^{m}$.
\end{corollary}

The next result describe this correspondence in more detail.

\begin{theorem}
Let $K$ be an incarnating set for a pair $Y\hookrightarrow l_{1}$Vectors $%
\left( y_{i}\right) $ of $K$ are parallel to ribs of the unit ball $%
B(Y^{\ast })$. A length of a sum of all vectors that are parallel to a given
rib $r$ of $B(Y^{\ast })$ is equal to the length of the rib $r$.
\end{theorem}

\begin{proof}
This result may be obtained as a consequence of [6]. Below will be presented
its direct proof by the induction on a dimension.

If $dim(K^{\prime })=1$ then $K^{\prime }=[-y,y]$ for some $y\in \mathbb{R}$%
; $K=\{-y,y\}$ and the theorem is true trivially.

Assume that the theorem is true for $dim(K\prime )=n$.

Let $K\in \mathbb{R}^{n+1}$; $dim(K^{\prime })=n+1$.

An $n$-dimensional subspace $E\hookrightarrow \mathbb{R}^{n+1}$ is said to
be a support subspace (for $K^{\prime }$) if $E\cap K$ is complete in $E$ \
A set of all support subspaces will be denoted by $E(K)$.

Let $E\in E(K)$; $r_{E}\in \mathbb{R}^{n+1}$; $r_{E}\neq 0$. Let $r_{E}$ be
orthogonal to $E$. Let
\begin{equation*}
x_{E}=\sum \{y_{i}\in K:\left\langle y_{i},r_{E}\right\rangle >0\}.
\end{equation*}

Then
\begin{equation*}
V_{E}=\left( E\cap K\right) ^{\prime }+x_{E}=\{z+x_{E}:z\in \left( E\cap
K\right) ^{\prime }\}
\end{equation*}
is a facet of $K^{\prime }$.

Indeed, $dim(K^{\prime })=n+1$; $dim(V_{E})=n$; $V_{E}\subset K$. At the
same time $E+x_{E}$ is an affine manifold of the minimal dimension that
contains $V_{E}$. By the definition of $K^{\prime }$, $\left( E+x_{E}\right)
\cap K^{\prime }=V_{E}$. Hence, $V_{E}$ is one of maximal $n$-dimensional
subsets of $K^{\prime }$. Moreover,
\begin{equation*}
V_{E}=\{z\in K^{\prime }:\left\langle z,r_{E}\right\rangle =\left\langle
x_{E},r_{E}\right\rangle =\underset{x\in K^{\prime }}{\sup }\left\langle
x,r_{E}\right\rangle \}.
\end{equation*}
So, to any support subspace $E\in E(K)$ corresponds a facet $V_{E}=\left(
E\cap K\right) ^{\prime }+x_{E}$.

Conversely, let $V$ be a facet of $K^{\prime }$; $r_{V}$ be an external
perpendicular to $V$. Let
\begin{eqnarray*}
x_{V} &=&\sum \{y\in K:\left\langle y,r_{V}\right\rangle >0\}; \\
E_{V} &=&\{x\in \mathbb{R}^{n+1}:\left\langle x,r_{V}\right\rangle =0\}.
\end{eqnarray*}
Then $V=\left( E_{V}\cap K\right) ^{\prime }+x_{V}$ and $E_{V}$ is a support
subspace.

Let $E\in E(K).$By the induction supposition, $E\cap K$ may be reconstructed
by the facet $V_{E}=\left( E\cap K\right) ^{\prime }+x_{V}$.

Since $K=\cup \{E\cap K:E\in E(K)\}$, then $K$ also may be reconstructed by
its facets. Notice that every rib belongs to some facet of $K$.
\end{proof}

No we are ready to prove the amalgamation property for $\frak{M}(\left(
l_{1}\right) ^{f})$.

\begin{theorem}
The set $\frak{M}(\left( l_{1}\right) ^{f})$ has the amalgamation property.
\end{theorem}

\begin{proof}
We show that $H(l_{1})$ has the amalgamation property. Obviously, a closure $%
\overline{H(l_{1})}=\frak{M}(\left( l_{1}\right) ^{f})$ also has this
property.

Let $Y\hookrightarrow l_{1}$ and $Z\hookrightarrow l_{1}$. According to
preceding results these pairs determine corresponding incarnating sets $%
K_{Y}=\left( y_{i}\right) $ and $K_{Z}=\left( z_{j}\right) $. Let $X$ be
isometricaly embedded into $Y$ (by operator $i:X\rightarrow Y)$ and into $Z$
(by operator $j:X\rightarrow Z$).

Consider a conjugate $X^{\ast }$ and for every rib of its unit ball $%
B(X^{\ast })$ choose subsets $\left( y_{i}^{r}\right) \subset K_{Y}=\left(
y_{i}\right) $ and $\left( z_{j}^{r}\right) \subset K_{Z}=\left(
z_{j}\right) $ that contains vectors parallel to $r$. Consider a set $%
K_{XYZ}=\left( x_{ij}\right) =\{c_{ij}^{r}\}$, where $c_{ij}^{r}=y_{i}\left%
\| z_{j}\right\| /\left\| r\right\| $; $\left\| r\right\| $ denotes a length
of a rib $r$. Clearly, $K_{XYZ}$\ is an incarnating set for a subspace $W$
of $l_{1}$ that contains isometric images of both $Y$ and $Z$ (say, $%
Y^{\prime }$ and $Z^{\prime }$) such that their intersection $y^{\prime
}\cap Z^{\prime }=X^{\prime }$ is exactly isometric to $X$; moreover,
corresponding isometries $u_{Y}:Y\rightarrow Y^{\prime }$ and $%
u_{Z}:Z\rightarrow Z^{\prime }$ transform embeddings $i:X\rightarrow Y$ and $%
j:X\rightarrow Z$ in identical embeddings of $X^{\prime }$ into $Y^{\prime }$
and into $Z^{\prime }$ respectively. Clearly this prove the theorem.
\end{proof}

\begin{remark}
Apparently, the same result may be obtained in another way by using one W.
Rudin's result [7] (see also [8]). It may be noted that Rudin examined the
complex case. Analogous result for real case is contained in [9].
\end{remark}

\section{Quotient closed and divisible classes of finite equivalence}

In this section it will be shown how to enlarge a Minkowski's base $\frak{M}%
(X^{f})$ of a certain $B$-convex class $X^{f}$ to obtain a set $\frak{N}$,
which will be a Minkowski's base $\frak{M}(W^{f})$ for some class $W^{f}$
that holds the $B$-convexity (and, also, $l_{p}$-spectrum) of the
corresponding class $X^{f}$ and has the amalgamation property.

\begin{definition}
A class $X^{f}$ (and its Minkowski's base $\frak{M}(X^{f})$) is said to be
divisible if some (equivalently, any) space $Z\in X^{f}$ is finitely
representable in any its subspace of finite codimension.
\end{definition}

\begin{definition}
Let $\{X_{i}:i\in I\}$ be a collection of Banach spaces. A space
\begin{equation*}
l_{2}\left( X_{i},I\right) =\left( \sum \oplus \{X_{i}:i\in I\}\right) _{2}
\end{equation*}
is a Banach space of all families $\{x_{i}\in X_{i}:i\in I\}=\frak{x}$, with
a finite norm
\begin{equation*}
\left\| \frak{x}\right\| _{2}=\sup \{(\sum \{\left\| x_{i}\right\|
_{X_{i}}^{2}:i\in I_{0}\})^{1/2}:I_{0}\subset I;\text{ }card\left(
I_{0}\right) <\infty \}.
\end{equation*}
\end{definition}

\begin{example}
Any Banach space $X$ may be isometricaly embedded into a space
\begin{equation*}
l_{2}(X)=(\sum\nolimits_{i<\infty }\oplus X_{i})_{2},
\end{equation*}
where all $X_{i}$'s are isometric to $X$. Immediately, $l_{2}(X)$ generates
a divisible class $\mathsf{D}_{2}(X^{f})=\left( l_{2}(X)\right) ^{f}$ which
with the same $l_{p}$-spectrum as $X^{f}$ and is superreflexive if and only
if $X^{f}$ is superreflexive too.
\end{example}

\begin{remark}
The procedure $\mathsf{D}_{2}:X^{f}\rightarrow \left( l_{2}(X)\right) ^{f}$
may be regarded as a closure operator on the partially ordered set $f\left(
\mathcal{B}\right) $. Indeed, it is

\begin{itemize}
\item  Monotone, i.e., $X^{f}<_{f}\mathsf{D}_{2}(X^{f})$;

\item  Idempotent, i.e., $\mathsf{D}_{2}(X^{f})=\mathsf{D}_{2}(\mathsf{D}%
_{2}(X^{f}))$;

\item  Preserve the order: $X^{f}<_{f}Y^{f}$ $\Longrightarrow \mathsf{D}%
_{2}(X^{f})<_{f}\mathsf{D}_{2}(Y^{f})$.
\end{itemize}

It is of interest that extreme points of $f\left( \mathcal{B}\right) $ are
stable under this procedure:
\begin{equation*}
\mathsf{D}_{2}(\left( l_{2}\right) ^{f})=\left( l_{2}\right) ^{f};\text{ }%
\mathsf{D}_{2}(\left( c_{0}\right) ^{f})=\left( c_{0}\right) ^{f}.
\end{equation*}
\end{remark}

To distinguish between general divisible classes and classes of type $%
\mathsf{D}_{2}(X^{f})$, the last ones will be called \textit{2-divisible
classes}.

\begin{definition}
A class $X^{f}$ (and its Minkowski's base $\frak{M}(X^{f})$) is said to be
quotient closed if for any $A\in \frak{M}(X^{f})$ and its subspace $%
B\hookrightarrow A$ the quotient $A/B$ belongs to $\frak{M}(X^{f})$.
\end{definition}

Let $K\subseteq \frak{M}$ be a class of finite dimensional Banach spaces
(recall, that isometric spaces are identified). Define operations $H$, $Q$
and $\ast $ that transform a class $K$ in another class - $H(K)$; $Q(K)$ or $%
(K)^{\ast }$ respectively. Namely, let

\begin{align*}
H(K) & =\{A\in\frak{M}:A\hookrightarrow B;\text{ \ }B\in K\} \\
Q(K) & =\{A\in\frak{M}:A=B/F;\text{ \ }F\hookrightarrow B;\text{ \ }B\in K\}
\\
(K)^{\ast} & =\{A^{\ast}\in\frak{M}:A\in K\}
\end{align*}

In words, $H(K)$ consists of all subspaces of spaces from $K$; $Q(K)$
contains all quotient spaces of spaces of $K$; $(K)^{\ast}$ contains all
conjugates of spaces of $K$.

The following theorem lists properties of these operations. In iteration of
the operations parentheses may be omitted.

Thus, $K^{\ast\ast}\overset{def}{=}\left( (K)^{\ast}\right) ^{\ast}$; $HH(K)%
\overset{def}{=}H(H(K))$ and so on.

\begin{theorem}
Any set $K$ of finite dimensional Banach spaces has the following properties:

\begin{enumerate}
\item  $K^{\ast \ast }=K$; $HH(K)=H(K)$; $QQ(K)=Q(K)$;

\item  $K\subset H(K)$; $K\subset Q(K)$;

\item  If $K_{1}\subset K_{2}$ then $H(K_{1})\subset H(K_{2})$ and $%
Q(K_{1})\subset Q(K_{2})$;

\item  $\left( H(K)\right) ^{\ast }=Q(K^{\ast })$; $\left( Q(K)\right)
^{\ast }=H(K^{\ast })$;

\item  $HQ(HQ(K))=HQ(K)$; $QH(QH(K))=QH(K)$.
\end{enumerate}
\end{theorem}

\begin{proof}
1, 2 and 3 are obvious.

4. If $A\in Q(K)$ then $A=B/E$ for some $B\in K$ and its subspace $E$. Thus,
$A^{\ast }$ is isometric to a subspace of $B^{\ast }.$ Hence, $A^{\ast }\in
H(B^{\ast })$, i.e., $A^{\ast }\in H(K^{\ast })$. Since $A$ is arbitrary, $%
\left( Q(K)\right) ^{\ast }\subseteq H(K^{\ast })$. Analogously, if $B\in K$
and $A\in H(B)$ then $A^{\ast }$ may be identified with a quotient $B^{\ast
}/A^{\perp }$, where $A^{\perp }$ is the annihilator of $A$ in $B^{\ast }$:
\begin{equation*}
A^{\perp }=\{f\in B^{\ast }:\left( \forall a\in A\right) \text{ \ }\left(
f\left( a\right) =0\right) \}.
\end{equation*}

Hence, $A^{\ast }\in (Q(K^{\ast }))^{\ast }$, and thus $\left( H(K)\right)
^{\ast }\subseteq Q(K^{\ast })$.

From the other hand,
\begin{align*}
H(K^{\ast })& =\left( H(K^{\ast })\right) ^{\ast \ast }\subseteq \left(
Q(K^{\ast \ast })\right) ^{\ast }=\left( Q(K)\right) ^{\ast }; \\
Q(K^{\ast })& =\left( Q(K^{\ast })\right) ^{\ast \ast }\subseteq \left(
H(K^{\ast \ast })\right) ^{\ast }=\left( H(K)\right) ^{\ast }.
\end{align*}

5. Let $A\in HQ(K)$. Then $A$ is isometric to a subspace of some quotient
space $E/F$, where $E\in K$; $F\hookrightarrow E$. If $B$ is a subspace of $%
A $ then $\left( A/B\right) ^{\ast }=\left( E/F\right) ^{\ast }/B^{\perp }$,
i.e. $\left( Q\left( HQ(K)\right) \right) ^{\ast }\subseteq Q\left(
Q(K)^{\ast }\right) $. Because of
\begin{align*}
Q\left( HQ(K)\right) & =\left( Q\left( HQ(K)\right) \right) ^{\ast \ast
}\subseteq (Q\left( Q(K)^{\ast }\right) )^{\ast } \\
& \subseteq H(Q(K)^{\ast \ast })=HQ(K),
\end{align*}

we have
\begin{equation*}
H(Q(H(Q(K))))\subseteq H(H(Q(K)))=HQ(K).
\end{equation*}

Analogously, if $A\in QH(K)$, then $A$ is isometric to a quotient space $F/E$%
, where $F\in H(B)$ for some $B\in K$ and $E\hookrightarrow B$. If $W\in
H(A) $, i.e., if $W\in H(F/E)$, then $W^{\ast }=(F/E)^{\ast }/W^{\perp }$
and $(F/E)^{\ast }$ is isometric to a subspace $E^{\perp }$ of $F^{\ast }\in
(H(B))^{\ast }$. Thus, $(H(QH(K)))^{\ast }\subseteq H((H(K))^{\ast })$ and
\begin{align*}
H\left( QH(K)\right) & =\left( H\left( QH(K)\right) \right) ^{\ast \ast
}\subseteq (H\left( H(K)^{\ast }\right) )^{\ast } \\
& \subseteq Q(H(K)^{\ast \ast })=QH(K).
\end{align*}

Hence,
\begin{equation*}
Q(H(Q(HQ(K))))\subseteq Q(Q(H(K)))=QH(K).
\end{equation*}

Converse inclusion follows from 2.
\end{proof}

Consider for a given 2-divisible class $X^{f}$ its Minkowski's base $\frak{M}%
(X^{f})$ and enlarge it by addition to $\frak{M}(X^{f})$ of all quotient
spaces of spaces from $\frak{M}(X^{f})$ and all their subspaces. In the
formal language, consider a set $H(Q(\frak{M}(X^{f})))$.

For any class $W^{f}$ the set $\frak{N}=\frak{M}(W^{f})$ has following
properties:

(\textbf{C}) $\frak{N}\ $is a closed subset of the Minkowski's space $\frak{M%
}$;

(\textbf{H}) If $A\in\frak{N}$ and $B\in H(A)$ then $B\in\frak{N}$;

(\textbf{A}$_{0}$) For any $A$, $B\in\frak{N}$ there exists $C\in\frak{N}$
such that $A\in H(C)$ and $B\in H(C)$.

Clearly, $H(Q(\frak{M}(X^{f})))$ has properties (\textbf{C}), (\textbf{H})
and (\textbf{A}$_{0}$).

\begin{theorem}
Let $\frak{N}$ be a set of finite dimensional Banach spaces; $\frak{N}%
\subset \frak{M}$. If $\frak{N}$ has properties (\textbf{C}), (\textbf{H})
and (\textbf{A}$_{0}$) then there exists a class $X^{f}$ such that $\frak{N}=%
\frak{M}(X^{f})$.
\end{theorem}

\begin{proof}
Since (\textbf{A}$_{0}$), all spaces from $\frak{N}$ may be directed in an
inductive isometric system. Its direct limit $W$ generates a class $W^{f}$
with $\frak{M}(W^{f})\supseteq \frak{N}$. As it follows from the property (%
\textbf{C}), $\frak{M}(W^{f})$ contains no spaces besides that of $\frak{N}$%
. Hence, $\frak{M}(W^{f})=\frak{N}$.

Another proof may be given by using an ultraproduct of all spaces from $%
\frak{N}$ by the ultrafilter, coordinated with a partial ordering on $\frak{N%
}$, which is generated by the property (A$_{0}$).
\end{proof}

Let $X$ be a $B$-convex Banach space.

Let $X=l_{2}(Y)$ (and, hence, $X^{f}=\mathsf{D}_{2}(Y^{f})$). Consider the
Minkowski's base $\frak{M}(X^{f})$ and its enlargement $H(Q(\frak{M}%
(X^{f})))=HQ\frak{M}(X^{f})$.

\begin{theorem}
There exists a Banach space $W$ such that $HQ\frak{M}(X^{f})=\frak{M}(W^{f})$%
.
\end{theorem}

\begin{proof}
Obviously, $HQ\frak{M}(X^{f})$ has properties (H) and (C).

Since $\frak{M}(X^{f})$ is 2-divisible, then for any $A$,$\ B\in \frak{M}%
(X^{f})$ a space $A\oplus _{2}B$ belongs to $\frak{M}(X^{f})$ and, hence, to
$HQ\frak{M}(X^{f})$.

If $A,B\in Q\frak{M}(X^{f})$ then $A=F/F_{1}$; $B=E/E_{1}$ for some $E,F\in Q%
\frak{M}(X^{f})$.

Clearly, $F/F_{1}\oplus _{2}E$ is isometric to a space ($F\oplus
_{2}E)/F_{1}^{\prime }$, where
\begin{equation*}
F_{1}^{\prime }=\{(f,0)\in F\oplus E:f\in F_{1}\}
\end{equation*}
and, hence, belongs to $Q\frak{M}(X^{f})$.

Thus, $F/F_{1}\oplus _{2}E/E_{1}$ =($F/F_{1}\oplus _{2}E)/E_{1}^{\prime }$,
where
\begin{equation*}
E_{1}^{\prime }=\{(o,e)\in F/F_{1}\oplus _{2}E:e\in E_{1}\},
\end{equation*}
and, hence, belongs to $Q\frak{M}(X^{f})$ too.

If $A,B\in HQ\frak{M}(X^{f})$ then $A\hookrightarrow E$, $B\hookrightarrow F$
for some $E,F\in Q\frak{M}(X^{f})$.

Since $E\oplus _{2}F\in Q\frak{M}(X^{f})$ also $A\oplus _{2}B\in HQ\frak{M}%
(X^{f})$. Thus, $HQ\frak{M}(X^{f})$ has a property (A$_{0}$). A desired
result follows from the preceding theorem.
\end{proof}

\begin{definition}
Let $X$ be a Banach space, which generates a class $X^{f}$\ of finite
equivalence. A class $\ast \ast (X^{f})$ is defined to be a such class $%
W^{f} $ that $\frak{M}(W^{f})=HQ\frak{M}(\mathsf{D}_{2}Y^{f})$
\end{definition}

Clearly, $W^{f}$ is quotient closed. Obviously, $X^{f}<_{f}W^{f}$. It will
be said that $W^{f}$ is a result of a procedure $\ast \ast $ that acts on $f(%
\mathcal{B})$.

\begin{remark}
The procedure $\ast \ast :X^{f}\rightarrow \ast \ast \left( X^{f}\right)
=W^{f}$ may be regarded as a closure operator on the partially ordered set $%
f\left( \mathcal{B}\right) $. Indeed, it is

\begin{itemize}
\item  Monotone, i.e., $X^{f}<_{f}\ast \ast (X^{f})$;

\item  Idempotent, i.e., $X^{f}=\ast \ast (\ast \ast (X^{f}))$;

\item  Preserve the order: $X^{f}<_{f}Y^{f}$ $\Longrightarrow \ast \ast
(X^{f})<_{f}\ast \ast (Y^{f})$.
\end{itemize}

It is of interest that extreme points of $f\left( \mathcal{B}\right) $ are
stable under this procedure:
\begin{equation*}
\ast \ast (\left( l_{2}\right) ^{f})=\left( l_{2}\right) ^{f};\text{\ }\ast
\ast (\left( c_{0}\right) ^{f})=\left( c_{0}\right) ^{f}.
\end{equation*}
\end{remark}

\begin{theorem}
For any Banach space $X$ a class $\ast \ast \left( X^{f}\right) $ is
2-divisible.
\end{theorem}

\begin{proof}
Let $\frak{N}=\frak{M}(\ast \ast (X^{f}))$\ . Since for any pair $A,B\in
\frak{N}$ \ their $l_{2}$-sum belongs to $\frak{N}$, then by an induction, $%
(\sum\nolimits_{i\in I}\oplus A_{i})_{2}\in \frak{N}$ for any finite subset $%
\{A_{i}:i\in I\}\subset \frak{N}$.

Hence, any infinite direct $l_{2}$-sum $(\sum\nolimits_{i\in I}\oplus
A_{i})_{2}$ is finite representable in $\ast \ast \left( X^{f}\right) $. Let
$\{A_{i}:i<\infty \}\subset \frak{N}$ \ is dense in $\frak{N}$. Let $%
Y_{1}=(\sum\nolimits_{i<\infty }\oplus A_{i})_{2}$; $Y_{n+1}=Y_{n}\oplus
_{2}Y_{1}$; $Y_{\infty }=\underset{\rightarrow }{\lim }Y_{n}$. Clearly, $%
Y_{\infty }=l_{2}\left( Y\right) $ and belongs to $\ast \ast \left(
X^{f}\right) $.
\end{proof}

Let $\star :f\left( \mathcal{B}\right) \rightarrow f\left( \mathcal{B}%
\right) $ be one else procedure that will be given by following steps.

Let $X\in \mathcal{B}$; $Y^{f}=\mathsf{D}_{2}\left( X^{f}\right) $. Let $%
\frak{Y}_{0}$ be a countable dense subset of $\frak{M}(Y^{f})$; $\frak{Y}%
_{0}=\left( Y_{n}\right) _{n<\infty }$. Consider a space $Z=(\sum_{n<\infty
}\oplus Y_{n})_{2}$ and its conjugate $Z^{\ast }$. $Z^{\ast }$ generates a
class $\left( Z^{\ast }\right) ^{f}$ which will be regarded as a result of
acting of a procedure $\star :X^{f}\rightarrow \left( Z^{\ast }\right) ^{f}$%
. Since $\left( Z^{\ast }\right) ^{f}$ is 2-divisible, iterations of the
procedure $\star $ are given by following steps: Let $\frak{Z}_{0}$ is a
countable dense subset of $\frak{M}(\left( Z^{\ast }\right) ^{f})$; $\frak{Z}%
_{0}=\left( Z_{n}\right) _{n<\infty }$. Consider a space $W=(\sum_{n<\infty
}\oplus Z_{n})_{2}$ and its conjugate $W^{\ast }$. Clearly,\ $W^{\ast }$
generates a class
\begin{equation*}
\left( W^{\ast }\right) ^{f}=\star \left( Z^{\ast }\right) ^{f}=\star \star
\left( X^{f}\right) .
\end{equation*}

\begin{theorem}
For any Banach space $X$ classes $\ast \ast \left( X^{f}\right) $ and $\star
\star \left( X^{f}\right) $ are identical.
\end{theorem}

\begin{proof}
From the construction follows that $H(Z^{\ast })=\left( QH(l_{2}(X)\right)
^{\ast }$ and that
\begin{align*}
H(W^{\ast })& =(QH(Z^{\ast }))^{\ast }=(Q(QH(l_{2}(X)))^{\ast })^{\ast } \\
& =H(QH(l_{2}(X)))^{\ast \ast }=HQH(l_{2}(X)).
\end{align*}

Hence,
\begin{equation*}
\frak{M}(\left( W^{\ast }\right) ^{f})=HQ\frak{M}(Y^{f})=HQ\frak{M}(\mathsf{D%
}_{2}(X^{f})).
\end{equation*}
\end{proof}

\begin{theorem}
Let $X$ be a Banach space, which generates a class of finite equivalence $%
X^{f}$. If $X$ is B-convex, then the procedure $\ast \ast $ maps $X^{f}$ to
a class $\ast \ast \left( X^{f}\right) $ with the same lower and upper
bounds of its $l_{p}$-spectrum as $X^{f}$. If $X$ is superreflexive then $%
\ast \ast \left( X^{f}\right) $ is superreflexive too.
\end{theorem}

\begin{proof}
Obviously, if $p\in S(X)$ and $p\leq 2$ then the whole interval $[p,2]$
(that may consist of one point) belongs to $S(X)$ and, hence, to $S(\ast
\ast \left( X^{f}\right) )$. If $p\in S(X)$ and $p>2$ then, by duality, $%
p/(p-1)\in S\left( \star \left( X^{f}\right) \right) $; hence $%
[p/(p-1),2]\subset S\left( \star \left( X^{f}\right) \right) $ and, thus $%
[2,p]\subset S(X)\subset S(\ast \ast \left( X^{f}\right) )$. If $p\notin
S(X) $ then $p/(p-1)\notin S\left( \star \left( X^{f}\right) \right) $ by
its construction and, hence, $p\notin S(\ast \ast \left( X^{f}\right) )$ by
the same reason. Hence, if
\begin{equation*}
p_{X}=\inf S(X)\text{; \ }q_{X}=\sup S(X)
\end{equation*}
then
\begin{equation*}
S(\ast \ast \left( X^{f}\right) )=[p_{X},q_{X}]=[\inf S(X),\sup S(X)].
\end{equation*}

The second assertion of the theorem is also obvious.
\end{proof}

\begin{theorem}
Let $X$ be a Banach space, which generates a class of finite equivalence $%
X^{f}$. If $X$ is of type $p>1$ and of cotype $q<\infty $, then the
procedure $\ast \ast $ maps $X^{f}$ to a class of the same type and cotype.
\end{theorem}

\begin{proof}
According to [7], for any $B$-convex Banach space $X$ of nontrivial type $p$
(and, hence, cotype $q$), its conjugate $X^{\ast }$ is of type $p^{\prime
}=q/(q-1)$ and of cotype $q^{\prime }=p/(p-1)$ (and, hence, is $B$ convex
too). Obviously, $\star \left( X^{f}\right) $ is of type $p^{\prime
}=q/(q-1) $ and of cotype $q^{\prime }=p/(p-1)$ and thus $\star \star \left(
X^{f}\right) =\ast \ast \left( X^{f}\right) $ is of type $p$ and of cotype%
\textit{\ }$q$.
\end{proof}

\section{ The amalgamation property of classes $\ast \ast X^{f}$}

Classes of kind $\ast \ast X^{f}$ have a more powerful property then the
amalgamation one.

\begin{definition}
Let $X\in \mathcal{B}$ \ generates a class $X^{f}$\ with a Minkowski's base $%
\frak{M}(X^{f})$. It will be said that $\frak{M}(X^{f})$ (and the whole
class $X^{f}$) has the isomorphic amalgamation property if for any fifth $%
\left\langle A,B_{1},B_{2},i_{1},i_{2}\right\rangle $, where $A$, $B_{1}$, $%
B_{2}\in \frak{M}(X^{f})$; $i_{1}:A\rightarrow B_{1}$ and$\
i_{2}:A\rightarrow B_{2}$ are isomorphic embeddings, there exists a triple $%
\left\langle j_{1},j_{2},F\right\rangle $, where $F\in \frak{M}(X^{f})$ and$%
\ j_{1}:B_{1}\rightarrow F$; $j_{2}:B_{2}\rightarrow F$ are isometric
embedding, such that $j_{1}\circ i_{1}=j_{2}\circ i_{2}$.
\end{definition}

\begin{remark}
It is easy to show that $\frak{M}$ has the isomorphic amalgamation property.
\end{remark}

\begin{theorem}
Any 2-divisible quotient closed class $X^{f}$ has the isomorphic
amalgamation property.
\end{theorem}

\begin{proof}
Let $A$, $B_{1}$, $B_{2}\in \frak{M}(X^{f})$; $i_{1}:A\rightarrow B_{1}$ and
$i_{2}:A\rightarrow B_{2}$ are isomorphic embeddings. Since $X^{f}$ is
2-divisible, $C=B_{1}\oplus _{2}B_{2}\in \frak{M}(X^{f})$. Consider a
subspace $H$ of $C$ that is formed by elements of kind $\left(
i_{1}a,-i_{2}a\right) $, where $a$ runs $A$:
\begin{equation*}
H=\{\left( i_{1}a,-i_{2}a\right) :a\in A\}.
\end{equation*}

Consider a quotient $W=C/H$. Since $X^{f}$ is quotient closed, $W\in \frak{M}%
(X^{f})$. Let $h:C\rightarrow W$ be a standard quotient map. Let $%
j_{1}:B_{1}\rightarrow W$ and $j_{2}:B_{2}\rightarrow W$ are given by:
\begin{align*}
j_{1}(b_{1})& =h(b_{1},0)\text{ \ \ \ for }b_{1}\in B_{1}; \\
j_{2}(b_{2})& =h(0,b_{2})\text{ \ \ \ for }b_{2}\in B_{2}.
\end{align*}

It is clear that $j_{1}$ and $j_{2}$ are isometric embeddings such that $\
j_{1}\circ i_{1}=j_{2}\circ i_{2}$.
\end{proof}

\begin{remark}
It may be shown that any quotient closed divisible (not necessary
2-divisible) class $X^{f}$ enjoys the isomorphic amalgamation property too.
For aims of the article the preceding result is satisfactory.
\end{remark}

\begin{remark}
Certainly, the isomorphic amalgamation property implies the amalgamation
property for any set $\frak{M}(X^{f})$. These properties are not equivalent.
It may be shown that the class $\frak{M}(\left( l_{1}\right) ^{f})$ has no
the isomorphic amalgamation property. Reasons of this fact will be described
later.
\end{remark}

Recall a definition, which is due to V.I. Gurarii [7].

\begin{definition}
Let $X$ be a Banach space; $\mathcal{K}$ be a class of Banach spaces. $X$ is
said to be a space of almost universal disposition with respect to $\mathcal{%
K}$ if for every pair of spaces $A$, $B$ of $\mathcal{K}$, such that $A$ is
a subspace of $B$ ($A\hookrightarrow B$), every $\varepsilon >0$ and every
isomorphic embedding $i:A\rightarrow X$ there exists an isomorphic embedding
$\hat{\imath}:B\rightarrow X$, which extends $i$ (i.e., $\hat{\imath}|_{A}=i$%
) and such, then
\begin{equation*}
\left\| \hat{\imath}\right\| \left\| \hat{\imath}^{-1}\right\| \leq
(1+\varepsilon )\left\| i\right\| \left\| i^{-1}\right\| .
\end{equation*}
\end{definition}

In the construction of the classical Gurarii's space of almost universal
disposition with respect to $\frak{M}$ the main role plays that $\frak{M}$
has the isomorphic amalgamation property.

The proof of a following results is almost literally repeats the Gurarii's
proof [7] on existence of a space of almost universal disposition with
respect to $\frak{M}$. Only changes that need to be made are: a substitution
of a set $\frak{M}$ \ with $\frak{M}(X^{f})$ for a given class $X^{f}$ and
using instead the mentioned above isomorphic amalgamation property of a set $%
\frak{M}$ \ the same property of a set $\frak{M}(X^{f})$. For this reason
its proof will be omitted.

\begin{theorem}
Any $\ast \ast $-closed class $X^{f}$ (as an arbitrary class of finite
equivalence, whose Minkowski's set $\frak{M}(X^{f})$ has the isomorphic
amalgamation property) contains a separable space $E_{X}$ of almost
universal disposition with respect to a set $\frak{M}(X^{f})$.

This space is unique up to almost isometry and is almost isotropic (in an
equivalent terminology, has an almost transitive norm: for any two elements $%
a$, $b\in E_{X}$, such that $\left\| a\right\| =\left\| b\right\| $ and
every $\varepsilon >0$ there exists an automorphism $u=u(a,b,\varepsilon
):E_{X}\overset{onto}{\rightarrow }E_{X}$ such that $\left\| u\right\|
\left\| u^{-1}\right\| \leq 1+\varepsilon $ and $ua=b$).

This space is an approximative envelope of a class $X^{f}$: for any $%
\varepsilon >0$ every separable Banach space which is finitely representable
in $X^{f}$ is $(1+\varepsilon )$-isomorphic to a subspace of $E_{X}$.
\end{theorem}

If $\frak{M}(X^{f})$ does not enjoy the isomorphic amalgamation property but
has only the amalgamation property, the corresponding class $X^{f}$ contains
a space, which has properties, similar\ to that of the aforementioned\ space
$E_{X}$.

\begin{definition}
Let $X$ be a Banach space; $\mathcal{K}$ be a class of Banach spaces. $X$ is
said to be almost $\omega $-homogeneous with respect to $\mathcal{K}$ if for
any pair of spaces $A$, $B$ of $\mathcal{K}$ such that $A$ is a subspace of $%
B$ ($A\hookrightarrow B$), every $\varepsilon >0$ and every isometric
embedding $i:A\rightarrow X$ there exists an isomorphic embedding $\hat{%
\imath}:B\rightarrow X$, which extends $i$ (i.e., $\hat{\imath}|_{A}=i$) and
such, then
\begin{equation*}
\left\| \hat{\imath}\right\| \left\| \hat{\imath}^{-1}\right\| \leq
(1+\varepsilon ).
\end{equation*}
\end{definition}

If $X$ is almost $\omega $-homogeneous with respect to $H(X)$ it will be
simply referred as to an almost $\omega $-homogeneous space.

\begin{theorem}
Any class $X^{f}$ of finite equivalence, whose Minkowski's set $\frak{M}%
(X^{f})$ has the amalgamation property contains a separable almost $\omega $%
-homogeneous space $G_{X}$.

This space is unique up to almost isometry and is almost isotropic (in an
equivalent terminology, has an almost transitive norm).

This space is an approximative envelope of a class $X^{f}$: for any $%
\varepsilon >0$ every separable Banach space which is finitely representable
in $X^{f}$ is $(1+\varepsilon )$-isomorphic to a subspace of $E_{X}$.
\end{theorem}

\begin{proof}
The proof of the theorem also literally repeats the Gurarii's one. We
present a sketch of proof for a sake of completeness (below will be
presented another proof of this fact, based on a different idea).

Let us starting with any space $X_{0}$ of $X^{f}$. Consider a dense
countable subset $\{x_{i}:i<\infty ;$ $\left\| x_{i}\right\| =1\}\subset
X_{0}$. Any finite subset $N$ of $\mathbb{N}$ defines a finite dimensional
subspace $U_{N}=span\{x_{i}:i<\infty \}\hookrightarrow X_{0}$. Clearly there
are only countable number of spaces of kind $U_{N}$. Identifying isometric
subspaces from
\begin{equation*}
\{U_{N}:N\subset \mathbb{N};\text{ \ }card(N)<\infty \}
\end{equation*}
it will be obtained a dense countable subset $\left( Z_{i}\right) _{i<\infty
8}\subset H(X_{0})$. Consider a set $\mathcal{F}$ of all triples $%
\left\langle A,B,i\right\rangle $ where $A$, $B\in \frak{M}(X^{f})$ and
there exists an isometric embedding $i:A\rightarrow B$.

Let $n<m<\infty $ and\ let $\mathcal{F}\left( n,m\right) \subset \mathcal{F}$
be a subset of $\mathcal{F}$, which consists of such triples $\left\langle
A,B,i\right\rangle $ that $dim\left( A\right) =n$; $dim\left( B\right) =m$.

$\mathcal{F}\left( n,m\right) $ may be equipped with a metric
\begin{equation*}
\varrho \left( \left\langle A,B,i\right\rangle ,\left\langle
A_{1},B_{1},i_{1}\right\rangle \right) =\log d^{n,m}(\left\langle
A,B,i\right\rangle ,\left\langle A_{1},B_{1},i_{1}\right\rangle ),
\end{equation*}
where
\begin{equation*}
d^{n,m}(\left\langle A,B,i\right\rangle ,\left\langle
A_{1},B_{1},i_{1}\right\rangle )=\inf \{\left\| u\right\| \left\|
u^{-1}\right\| :u:B\rightarrow B_{1};\text{ \ }i\circ u=i_{1}\}
\end{equation*}
is a generalized Banach-Mazur distance. It is known [10] that $\left\langle
\mathcal{F}\left( n,m\right) ,\varrho \right\rangle $ is a compact metric
space. Hence, $\left\langle \mathcal{F},\varrho \right\rangle =\cup
_{n,m}\left\langle \mathcal{F}\left( n,m\right) ,\varrho \right\rangle $ is
a separable metric space and it may be chosen a countable dense subset $%
\mathcal{F}^{0}$ of $\mathcal{F}$. Without loss of generality it may be
assumed that all spaces that are presented in triples from $\mathcal{F}^{0}$%
are exactly those that belong to a defined before subset $\left(
Z_{i}\right) \subset H(X_{0})$.

Using the amalgamation property for $\frak{M}(X^{f})$ for any pair $%
A\hookrightarrow B\hookrightarrow X_{0}$ of subspaces of $X_{0}$ and any
isometric embedding $i:A\rightarrow C$, where $C\in \frak{M}(X^{f})$ it may
be constructed an extension of $X_{0}$, say, $X_{0}(A\hookrightarrow B;$ $%
i:A\rightarrow C)$, - a separable Banach space that contains $X_{0}$ and
contains a triple $(A,B,C^{\prime })$ where $A\hookrightarrow B$; $%
A\hookrightarrow C^{\prime }$ and a pair $A\hookrightarrow C^{\prime }$ is
isometric to the pair $iA\hookrightarrow C$ in a sense of the aforementioned
metric $\varrho $. It will be said that $X_{0}(A\hookrightarrow B;$ $%
i:A\rightarrow C)$ amalgamates over $X_{0}$ the $V$-formation $%
(A\hookrightarrow B;$ $i:A\rightarrow C)$.

Now, we proceed by an induction.

Let $\left( \frak{f}_{i}\right) _{i<\infty }$ be a numeration of all triples
from $\mathcal{F}^{0}$. Construct a sequence of Banach spaces $\left(
X_{i}\right) _{i<\infty }$ where will be presented the space $X_{0}$ itself
as a first step of the induction.

Let spaces $\left( X_{i}\right) _{i<\infty }$ \ be already constructed. As $%
X_{n+1}$ will be chosen the space $X_{n}(A\hookrightarrow B;$ $%
i:A\rightarrow C),$ where $A,B\in \left( U_{N}\right) _{N\subset \mathbb{N}}$%
; $A\hookrightarrow B\hookrightarrow X_{n}$ be exactly the $n$'s triple $%
\frak{f}_{n}$; $i:A\rightarrow C=\frak{f}_{m+1}$, where $m$ is the least
number of a triple from $\mathcal{F}^{0}$ that is of kind $i:A\rightarrow C$
(for a fixed $A$ and arbitrary $i$ and $C$) and is such that $X_{n}$ does
not amalgamate the $V$-formation $(A\hookrightarrow B;$ $i:A\rightarrow C)$.

Clearly, $X_{0}\hookrightarrow X_{1}\hookrightarrow X_{2}\hookrightarrow
...\hookrightarrow X_{n}\hookrightarrow ...$. Let $X_{\infty }^{(1)}=%
\overline{\cup X_{i}}=\underset{\rightarrow }{\lim }X_{i}$. Let us continue
the induction, started with $X_{\infty }^{(1)}$ instead of $X_{0}$. It will
be consequently constructed spaces $X_{\infty }^{(2)}\hookrightarrow
X_{\infty }^{(3)}\hookrightarrow ...\hookrightarrow X_{\infty
}^{(n)}\hookrightarrow ...$.

Their direct limit $\underset{\rightarrow }{\lim }X_{\infty }^{(n)}=G_{X}$
will be a desired space.

Proofs of the uniqueness of $G_{X}$ up to almost isometry and of the
property to be almost isotropic are the same as in [10].
\end{proof}

\section{Sub $B$-convex Banach spaces and classes}

A Banach space $X$ is said to be \textit{sub }$B$\textit{-convex} if it may
be represented as a direct sum $X=l_{1}\oplus _{1}W$, where $W$ is an
\textbf{infinite dimensional} $B$-convex space.

The $l_{1}$-sum is chosen for the convenience, all results remain true (up
to isomorphism) if instead $\oplus _{1}$ will be regarded any other norm $%
\mathcal{N}$ on the linear space $X\oplus l_{1}=\{(x,y):x\in X;$ $\ y\in Y\}$
(a corresponding Banach space may be denoted by $X\oplus _{\mathcal{N}}l_{1}$%
). Indeed, if the norm $\mathcal{N}$ is consistent with norms on $X$ and on $%
l_{1}$, i.e. if operators of canonical embeddings
\begin{eqnarray*}
u &:&X\rightarrow X\oplus Y;\text{ \ \ }ux=(x,0)\text{ \ for all }x\in X; \\
v &:&Y\rightarrow X\oplus Y;\text{ \ \ }uy=(0,y)\text{ \ for all }y\in l_{1}
\end{eqnarray*}
are continuous, then spaces $X\oplus _{1}l_{1}$ and $X\oplus _{\mathcal{N}%
}l_{1}$ are isomorphic.

Since $W$ is $B$-convex, $p_{W}=\inf S(W)>1$. It will be said that $\tau
(X)=p_{W}$ is a\textit{\ sub }$B$\textit{-type} of a sub $B$-convex space $%
X=l_{1}\oplus _{1}W$.

If a class $X^{f}$ contains a sub $B$-convex Banach space $Y$ then it also
will be referred as to the \textit{sub }$B$\textit{-convex class }$X^{f}$.
However, the value $\tau (Y)$ does not uniquely determined the sub $B$-type
of the class $X^{f}$. Indeed, if $l_{1}\oplus _{1}W$ and $p_{W}=p>1$ then,
since $l_{1}\oplus _{1}l_{r}<_{f}l_{1}$ for any $r\in (1,p)$, a space $%
l_{1}\oplus _{1}l_{r}\oplus _{1}W$ also belongs to the same class of finite
equivalence as $l_{1}\oplus _{1}W$; is also sub $B$-convex ($l_{r}\oplus
_{1}W=W_{1}$ obviously is $B$-convex) and $\tau \left( l_{1}\oplus
_{1}l_{r}\oplus _{1}W\right) =r<p_{W}=\tau (l_{1}\oplus _{1}W)$.

\begin{definition}
Let $X^{f}$ be a sub $B$-convex class. Then its sub $B$-type $\tau \left(
X^{f}\right) $ is given by
\begin{equation*}
\tau \left( X^{f}\right) =\sup \{\tau \left( Y\right) :Y\in X^{f}\text{ and }%
Y\text{ is sub }B\text{-convex\}.}
\end{equation*}
\end{definition}

\begin{remark}
Clearly, $\inf \{\tau \left( Y\right) :Y\in X^{f}$ and $Y$ is sub $B$-convex$%
\}=1$ for any sub $B$-convex class $X^{f}$. Notice that the space $l_{1}$
generates a sub $B$ convex class $\left( l_{1}\right) ^{f}$ (since $%
l_{1}\oplus _{1}l_{r}\in \left( l_{1}\right) ^{f}$ for any $1<r<2$).
Clearly, $\tau \left( \left( l_{1}\right) ^{f}\right) =2$. Notice also,
that, by the Dvoretzki theorem, $\tau \left( X^{f}\right) <2$ for every sub $%
B$-convex class $X^{f}$.
\end{remark}

More precise then by$\ \tau \left( X^{f}\right) $, sub $B$-convex class $%
X^{f}$ may be characterized by the notion of \textit{generative }(or \textit{%
determinative}) \textit{class} $\left( W_{X}\right) ^{f}$.

Clearly, if $l_{1}\oplus _{1}W\in X^{f}$ then $l_{1}\oplus _{1}W_{1}\in
X^{f} $ for any $W_{1}\sim _{f}W$. Moreover, if $l_{1}\oplus _{1}W\in X^{f}$
and $l_{1}\oplus _{1}U\in X^{f}$ and $V\in W^{f}\cap U^{f}$ then $%
l_{1}\oplus _{1}V\in X^{f}$ too.

\begin{definition}
A ($B$-convex) class $\left( W_{X}\right) ^{f}$ is said to be a
determinative class for a sub $B$-convex class $X^{f}$ provided
\begin{equation*}
\left( W_{X}\right) ^{f}=\cap \{W^{f}:l_{1}\oplus _{1}W\in X^{f}\}.
\end{equation*}
\end{definition}

Clearly, $\tau \left( X^{f}\right) =\tau \left( \left( W_{X}\right)
^{f}\right) $ and $\left( W_{X}\right) ^{f}$ is the least class of finite
equivalence (in a sense of partial order $<_{f}$) such that $l_{1}\oplus
_{1}W_{X}\in X^{f}$.

Notice that $\left( W_{X}\right) ^{f}\neq \varnothing $ since $\left(
l_{2}\right) ^{f}\in \cap _{i}\left( Z_{i}\right) ^{f}$ for any collections
of classes $\left( Z_{i}\right) ^{f}$ that are generated by infinite
dimensional spaces.

Since all $B$-convex classes $W^{f}$ are partially ordered by the relation $%
<_{f}$, it will be assumed that sub $B$-convex classes $X^{f}$ are also
ordered in a such way: $X^{f}\ll Y^{f}$ if $\left( W_{X}\right)
^{f}<_{f}\left( W_{X}\right) ^{f}$. Obviously, this order is identical to
the standard one: $X^{f}\ll Y^{f}$ if and only if $X^{f}<_{f}Y^{f}$.

The main goal of this section is to show how to enlarge a sub $B$-convex
class $X^{f}$ to a some class of finite equivalence $\left( G_{X}\right)
^{f} $\ (more precise, how to include the Minkowski's base $\frak{M}\left(
X^{f}\right) $ of a given sub $B$-convex class $X^{f}$ in a possibly minimal
set $\frak{M}\left( \left( G_{X}\right) ^{f}\right) $) with the amalgamation
property.

\begin{theorem}
Let $X^{f}$ be a sub $B$-convex class of finite equivalence. Then there
exists a class $Y^{f}$ with the following properties:

\begin{enumerate}
\item  $X^{f}<_{f}Y^{f}$

\item  $\frak{M}\left( Y^{f}\right) $ has the amalgamation property

\item  $Y^{f}$ is of the same cotype as $X^{f}$.
\end{enumerate}
\end{theorem}

\begin{proof}
It may be assumed that $X$ itself is sub $B$-convex, i.e. $X=l_{1}\oplus
_{1}W$, where $W$ is $B$-convex. Consider a class $\ast \ast \left(
W^{f}\right) $ and choose some $W_{1}\in \ast \ast \left( W^{f}\right) $.

Clearly, $X<_{f}l_{1}\oplus _{1}W_{1}=Z$ and $X$ is of the same cotype as $Z$%
.

To obtain the needed enlargement of $\frak{M}\left( X^{f}\right) $, let us
proceed by induction.

At the first step let $\frak{N}_{0}=\frak{M}\left( Z^{f}\right) $.
Certainly, $\frak{M}\left( X^{f}\right) \subseteq \frak{M}\left(
Z^{f}\right) $.

Assume that $\frak{N}_{0}$, $\frak{N}_{1}$, ... , $\frak{N}_{n}$ are already
defined. Define $\frak{N}_{n+1}^{0}$ to be an enlargement of $\frak{N}_{n}$
by addition to it all spaces of kind
\begin{equation*}
F=(B\oplus _{1}C)/N,
\end{equation*}
where $B\in \frak{M}\left( \left( l_{1}\right) ^{f}\right) $; $C\in \frak{N}%
_{n}\backslash \frak{M}\left( \left( l_{1}\right) ^{f}\right) $; $N=N\left(
A,u\right) $ is a subspace of $B\oplus _{1}C$ that depends on a given
subspace $A\hookrightarrow B$ and on isomorphic embedding $u:A\rightarrow C$%
, which is given by:
\begin{equation*}
N=N\left( A,u\right) =\{\left( a,-u(a)\right) :a\in A\}.
\end{equation*}

Let $\frak{N}_{n+1}=H\left( \frak{N}_{n+1}^{0}\right) $. Let
\begin{equation*}
\frak{N}_{\infty }=\cup \{\frak{N}_{n}:n<\infty \}.
\end{equation*}

Clearly, $\frak{N}_{\infty }$ (and, hence, its closure $\overline{\frak{N}%
_{\infty }}$) has the amalgamation property. To show this, proceed by
induction.

At the first step of induction consider a $V$-formation $v=\left\langle
A,B_{1},B_{2},i_{1},i_{2}\right\rangle $, whose spaces $A$, $B_{1}$, $B_{2}$
belongs to $\frak{N}_{0}$.

Certainly, $A=A^{\prime }\oplus _{1}A^{\prime \prime }$; $%
B_{1}=B_{1}^{\prime }\oplus _{1}B_{1}^{\prime \prime }$; $%
B_{2}=B_{2}^{\prime }\oplus _{1}B_{2}^{\prime \prime }$, where $A^{\prime }$%
; $B_{1}^{\prime }$; $B_{2}^{\prime }$ belongs to $\frak{M}\left( \left(
l_{1}\right) ^{f}\right) $; $A^{\prime \prime }$; $B_{1}^{\prime \prime }$; $%
B_{2}^{\prime \prime }\in \frak{M}\left( \ast \ast \left( Z\right)
^{f}\right) $.

If either $B_{1}^{\prime \prime }=B_{2}^{\prime \prime }=\{0\}$ or $%
B_{1}^{\prime }=B_{2}^{\prime }=\{0\}$ then $v$ is amalgamated in $\frak{N}%
_{0}$, since both $\frak{M}\left( \left( l_{1}\right) ^{f}\right) $ and $%
\frak{M}\left( \ast \ast \left( Z\right) ^{f}\right) $ has the amalgamation
property. In the opposite case, this $V$-formation has an amalgam in $\frak{N%
}_{1}$, as it follows from the definition of $\frak{N}_{1}$. On a general
case any $V$-formation $v$ in $\frak{N}_{\infty }$ has a root $A\in \frak{N}%
_{k}$ for some $k<\infty $ and spaces $B_{1}\in \frak{N}_{l}$, $B_{2}\in
\frak{N}_{m}$ respectively. If $v$ has no amalgam in $\frak{N}_{s}$, where $%
s=\max \{l,m\}$, then from the construction follows that $v$ has an amalgam
in $\frak{N}_{s+1}$. Indeed, in a such case, one of spaces $B_{1}$, $B_{2}$
belongs to $\frak{M}\left( \left( l_{1}\right) ^{f}\right) $ and other to $%
\frak{N}_{s}\backslash \frak{M}\left( \left( l_{1}\right) ^{f}\right) $.
Clearly, the space $F=(B_{1}\oplus _{1}B_{2})/N$, where $N=\{(i_{1}\left(
a\right) ,-i_{2}\left( a\right) :a\in A\}$, and isometric embeddings $%
j_{1}:B_{1}\rightarrow F$ and $j_{2}:B_{2}\rightarrow F$, which are given
by:
\begin{align*}
j_{1}(b_{1})& =h(b_{1},0)\text{ \ \ \ for }b_{1}\in B_{1}; \\
j_{2}(b_{2})& =h(0,b_{2})\text{ \ \ \ for }b_{2}\in B_{2},
\end{align*}
where $h:B_{1}\oplus _{1}B_{2}\rightarrow F$ is a standard quotient map,
form a triple $\left\langle j_{1},j_{2},F\right\rangle $, which is an
amalgam of $v$ because of $j_{1}\circ i_{1}=j_{2}\circ i_{2}$.

Clearly, the amalgamation property of $\frak{N}_{\infty }$ implies the
property (A$_{0}$). Thus, $\frak{N}_{\infty }$ has properties (H) and (A$%
_{0} $) and $\overline{\frak{N}_{\infty }}$ has properties (C), (H) and (A$%
_{0}$). Hence, there exists a class $Y^{f}$ whose Minkowski's base is
exactly $\overline{\frak{N}_{\infty }}$.

To show that $Y^{f}$ is of the same cotype as $X^{f}$ it is sufficient to
notice that for all spaces of kind $N=N\left( A,u\right) =\{\left(
a,-u(a)\right) :a\in A\}$ (that were used in the construction) their $p$%
-type constants $T_{p}\left( N\right) $, where $p=\tau (W)$, are uniformly
bounded. Assume that $X$ is of cotype $q$. According to J. Pisier's results
[11] (see also [1]) the $q$-cotype constants $C_{q}(A)$, where $A\in \frak{N}%
_{\infty }$ are also uniformly bounded.

Since the definition of cotype depends only on finite extracts of vectors, $%
Y^{f}$ is of cotype $q$.
\end{proof}

\begin{remark}
The class $Y^{f}$ that is constructed above depends on $W$. It will be said
that $Y^{f}$ is obtained from $X^{f}$ as a result of a procedure $\frak{Q}%
_{W}$, which is defined on a set of all sub $B$-convex classes:
\begin{equation*}
Y^{f}\overset{def}{=}\frak{Q}_{W}X^{f}.
\end{equation*}
If $Z=l_{1}\oplus _{1}U$ belongs to $X^{f}$ and $\tau (W)\neq \tau (U)$ then
classes $\frak{Q}_{W}X^{f}\ $and $\frak{Q}_{U}X^{f}$ are different, as will
be shown below.

If for $W$ will be chosen any space from the \textit{determinative} \textit{%
class} $\left( W_{X}\right) ^{f}$, say, $W_{X}$, the procedure $\frak{Q}%
_{W_{X}}X^{f}$ will be denoted simply by $\frak{Q}X^{f}$

If a class $Y^{f}$ is equal to $\frak{Q}_{W}X^{f}$ for some sub $B$-convex
class $X^{f}$ and $B$-convex $W$ it will be called a $\frak{Q}$-closed class.
\end{remark}

\begin{theorem}
Any $\frak{Q}$-closed class $X^{f}$ contains a separable almost $\omega $%
-homogeneous space $G_{X}$ which is unique up to almost isometry, is an
approximative envelope of a class $X^{f}$ and is almost isotropic.
\end{theorem}

\begin{proof}
Immediately follows from the preceding theorem and the theorem 15.
\end{proof}

\section{Some properties of classes $\frak{Q}_{W}X^{f}$}

Let $X^{f}$ be a sub $B$-convex class; $X=l_{1}\oplus _{1}W$; $E_{W}\in \ast
\ast (W^{f})$ be a space of almost universal disposition; $G_{X}$ be a
separable almost $\omega $-homogeneous space of $\frak{Q}_{W}X^{f}$. It may
be assumed that $E_{W}$ is isometric to a subspace of $G_{X}$. It will be
convenient to regard $E_{W}$ as a subspace of $G_{X}$.

\begin{theorem}
For any $A,$ $B\in \frak{M}\left( \frak{Q}_{W}X^{f}\right) $ and every
isometric embedding $i:A\rightarrow B$ each (linear) mapping $h:A\rightarrow
E_{W}$ may be extended to a mapping $\overline{h}:B\rightarrow G_{X}$ such
that $\overline{h}\circ i=h$ and $\left\| \overline{h}\right\| \leq c\left\|
h\right\| $ where the constant $c$ depends only on $W$.
\end{theorem}

\begin{proof}
Immediately follows from constructions of $\ast \ast (W^{f})$; $E_{W}$; $%
\frak{Q}_{W}X^{f}$ and $G_{X}$.
\end{proof}

\begin{definition}
(Cf. [12]). Let $X$, $Y$ be Banach spaces; $Y\hookrightarrow X$. $Y$ is said
to be a reflecting subspace of $X$, shortly: $Y\prec _{u}X$, if for every $%
\varepsilon >0$ and every finite dimensional subspace $A\hookrightarrow X$
there exists an isomorphic embedding $u:A\rightarrow Y$ such that $\left\|
u\right\| \left\| u^{-1}\right\| \leq 1+\varepsilon $, which is identical on
the intersection $Y\cap A$:
\begin{equation*}
u\mid _{Y\cap A}=Id_{Y\cap A}.
\end{equation*}
\end{definition}

Clearly, $Y\prec _{u}X$ implies that $Y\sim _{f}X$.

\begin{definition}
(Cf. [13]). A Banach space $E$ is said to be existentialy closed in a class $%
X^{f}$ if for any isometric embedding $i:E\rightarrow Z$ into an arbitrary
space $Z\in X^{f}$ its image $iE$ is a reflecting subspace of $Z$: $iY\prec
_{u}Z$.
\end{definition}

A class of all spaces $E$ that are existentialy closed in $X^{f}$ is denoted
by $\mathcal{E}\left( X^{f}\right) $. In [11] was shown that for any Banach
space $X$ the class $\mathcal{E}\left( X^{f}\right) $ is nonempty; moreover,
any $Y<_{f}X^{f}$ may be isometricaly embedded into some $E\in \mathcal{E}%
\left( X^{f}\right) $ of the dimension $dim(E)=\max \{dim(Y),\omega \}$.

\begin{theorem}
For any class $X^{f}$ such that $\frak{M}\left( X^{f}\right) $ has the
amalgamation property $E\in \mathcal{E}\left( X^{f}\right) $ if and only if $%
E\in X^{f}$ and $E$ is almost $\omega $-homogeneous.
\end{theorem}

\begin{proof}
Let $E\hookrightarrow Z\in X^{f}$ and $E$ be almost $\omega $-homogeneous.
Let $A\hookrightarrow Z$ be a finite dimensional subspace; $E\cap A=B$ and $%
\varepsilon >0$. Consider the identical embedding $id_{B}:B\rightarrow E$.
Since $B\hookrightarrow A$, $id_{B}$ may be extended to an embedding $%
u:A\rightarrow E$ with $\left\| u\right\| \left\| u^{-1}\right\| \leq
1+\varepsilon $. Thus, $E\in \mathcal{E}\left( X^{f}\right) $.

Conversely, let $E\in \mathcal{E}\left( X^{f}\right) $; $B\hookrightarrow E$%
; $B\hookrightarrow A$ and $A\in \frak{M}\left( X^{f}\right) $.

Consider a space $Z$ such that $E\hookrightarrow Z$; $B\hookrightarrow
A\hookrightarrow Z$. Such space exists because of the amalgamation property
of $\frak{M}\left( X^{f}\right) $. Since $E\in \mathcal{E}\left(
X^{f}\right) $, $E\prec _{u}Z$, i.e. there is an embedding $u:A\rightarrow E$
such that $\left\| u\right\| \left\| u^{-1}\right\| \leq 1+\varepsilon $,
which is identical on the intersection $E\cap A$: $u\mid _{E\cap
A}=Id_{E\cap A}$. Since $B\hookrightarrow A$ and $B\hookrightarrow E$, $%
B\hookrightarrow E\cap A$. Clearly, $u$ extends the identical embedding $%
id_{B}$. Since $A\hookrightarrow B$ and $\varepsilon $ are arbitrary, $E$ is
almost $\omega $-homogeneous.
\end{proof}

\begin{remark}
This theorem gives an alternative proof of the first part of theorem 15.
\end{remark}

\begin{theorem}
Let $X^{f}$ be a $\frak{Q}$-closed class; $G^{X}$ be the corresponding
separable almost $\omega $-homogeneous space. Then $G_{X}$ is not isomorphic
to a subspace of any conjugate separable Banach space.

Its conjugate $\left( G_{X}\right) ^{\ast }$ is nonseparable; moreover, $%
dim(\left( G_{X}\right) ^{\ast })>2^{\aleph _{0}}$.
\end{theorem}

\begin{proof}
Since $G_{X}$ is an approximative envelope of a class $X^{f},$it contains
(almost isometricaly) $L_{1}[0,1]$. Thus, the first part of theorem follows
from the classical result of I.M. Gelfand [14]. The second part is a
consequence of [15].
\end{proof}

\begin{remark}
It seems that the more direct way is to compute the Szlenk index of $G_{X}$.
Clearly, $G_{X}$ contains spaces $X_{n}=\left( \sum\nolimits_{k=1}^{n}\oplus
l_{2}\right) _{1}$ for any finite $n$ and $\left(
\sum\nolimits_{k=1}^{\infty }\oplus X_{k}\right) _{2}$ for any sequence $%
\left( X_{n}\right) $ of subspaces of $G_{X}$. Hence, its Szlenk index $%
Szl(G_{X})$ cannot be a countable ordinal and, according to [16] and [17], $%
G_{X}$ cannot be isomorphic to a subspace of any conjugate separable Banach
space, and its conjugate $\left( G_{X}\right) ^{\ast }$ is nonseparable
(apparently, the more powerful estimate for dimension of $\left(
G_{X}\right) ^{\ast }$ cannot be obtained on this way).

Notice, that in a similar way may be obtained the Gelfand theorem [14] and
some results of [18].
\end{remark}

Let $X,Y$ be Banach spaces. A Banach space of all (linear bounded) operators
$u:X\rightarrow Y$ will be denoted by $L(X,Y)$.

Let $p\in \lbrack 1,\infty )$; $u\in L(X,Y)$. $u$ is said to be

\begin{itemize}
\item  $p$-\textit{absolutely summing,} if there is a constant $\lambda >0$
such that
\begin{equation*}
(\sum\nolimits_{j<n}\left\| u\left( x_{i}\right) \right\| ^{p})^{1/p}\leq
\lambda (\sum\nolimits_{j<n}\left| \left\langle x_{i},f\right\rangle \right|
^{p})^{1/p}
\end{equation*}
for any $f\in X^{\ast }$ and any finite set $\{x_{i}:i<n;n<\infty \}\subset
X $.
\end{itemize}

Its $p$-\textit{absolutely summing norm }$\pi _{p}\left( u\right) $ is the
smallest constant $\lambda $.

A space of all $p$-absolutely summing operators from $X$ to $Y$ is denoted
by $AS_{p}\left( X,Y\right) $.

\begin{itemize}
\item  $p$-\textit{factorable,} if there exists a such measure $\mu $ and
operators $v\in B(X,L_{p}(\mu ))$; $w\in B(L_{p}(\mu ),Y^{\ast \ast })$ such
that $k_{Y}\circ u=w\circ v$, where $k_{Y}$ is the canonical embedding of $Y$
in its second conjugate $Y^{\ast \ast }$, and the symbol $\circ $ denotes
the composition of operators.
\end{itemize}

Its $\gamma _{p}$\textit{-norm }is given by
\begin{equation*}
\gamma _{p}\left( u\right) =\inf \{\left\| v\right\| \left\| w\right\|
:k_{Y}\circ u=w\circ v\}.
\end{equation*}

A space of all $p$-factorable operators from $X$ to $Y$ is denoted by $%
\Gamma _{p}\left( X,Y\right) $.

\begin{itemize}
\item  $p$-\textit{integral,} if there exists a such probability measure $%
\mu $ and such operators $v\in B(X,L_{\infty }(\mu ))$, $w\in B(L_{p}(\mu
),Y^{\ast \ast })$ that $k_{Y}\circ u=w\circ \varphi \circ v$, where $%
\varphi $ is an inclusion of $L_{\infty }(\mu )$ into $L_{p}(\mu )$.
\end{itemize}

Its $p$-\textit{integral norm} is given by
\begin{equation*}
\iota _{p}(u)=\inf \{\left\| v\right\| \left\| \varphi \right\| \left\|
w\right\| :k_{Y}\circ u=w\circ \varphi \circ v\}.
\end{equation*}

A space of all $p$-integral operators from $X$ to $Y$ is denoted by $%
I_{p}\left( X,Y\right) $.

\begin{itemize}
\item  \textit{Nuclear,} if there are sequences $(f_{n})_{n<\infty }\in
X^{\ast }$; $\left( y_{n}\right) _{n<\infty }\in Y$, such that for all $x\in
X$
\begin{equation*}
u\left( x\right) =\sum\nolimits_{n<\infty }f_{n}\left( x\right) y_{n}.
\end{equation*}
\end{itemize}

Its \textit{nuclear norm} $\nu _{1}(u)$ is given by
\begin{equation*}
\nu _{1}\left( u\right) =\inf \{\sum \left\| f_{n}\right\| \left\|
y_{n}\right\| :u\left( x\right) =\sum\nolimits_{n<\infty }f_{n}\left(
x\right) y_{n}\},
\end{equation*}
where the infimum is given over all possible representations of $u$.

A space of all nuclear operators from $X$ to $Y$ is denoted by $N_{1}\left(
X,Y\right) $.

\begin{theorem}
Let $X^{f}$ be a sub $B$-convex class that is generated by $X=l_{1}\oplus
_{1}W$; $G_{X}^{W}\ $be the almost $\omega $-homogeneous space for a class $%
\frak{Q}_{W}\left( (l_{1}\oplus _{1}W)^{f}\right) $. Let $\tau (W)=r$. Then
for every $p\in \lbrack r,2]$ there exists a such constant $\lambda $ that
for every Banach space $Y$ and any finite rank operator $u:Y\rightarrow
G_{X}^{W}$
\begin{equation*}
\iota _{p}\left( u\right) \leq \lambda \pi _{p}\left( u\right) .
\end{equation*}
\end{theorem}

\begin{proof}
Let $u:Y\rightarrow G_{X}^{W}$ be a finite rank operator and its $p$%
-absolutely summing norm is bounded by $c<\infty $: $\pi _{2}\left( u\right)
\leq c$. Consider the unit ball $S=B(Y^{\ast })$ of the space $Y^{\ast }$\
endowed with weak* topology.

Let $j:Y\rightarrow C(S)$ be the canonical embedding: $jy(y^{\prime
})=y^{\prime }(y)$ for $y\in Y$, $y^{\prime }\in S$.

Let $\mu $ be a measure on $S$; $i_{\mu ,p}:C(S)\rightarrow L_{p}\left(
S,\mu \right) $ be the natural embedding.

Notice that $\left( i_{\mu ,p}\circ j\right) Y=A_{0}$ is a linear subspace
of $L_{p}\left( S,\mu \right) $. Let $A$ be a closure of $A_{0}$ in the $%
L_{p}\left( S,\mu \right) $ metric. Let an operator $w_{0}:A_{0}\rightarrow
G_{X}^{W}$ be given by
\begin{equation*}
w_{0}\circ i_{\mu ,p}\circ j\left( z\right) =u(z);\text{ }z\in A_{0}.
\end{equation*}

According to the A. Pietsch theorem [20], the measure $\mu $ in a such way
that the operator $w_{0}$ will be continuous (in the $L_{p}\left( S,\mu
\right) $ metric) and, hence, may be extended to the operator $%
w:A\rightarrow G_{X}^{W}$.

Because of $A$ is finite dimensional, it may be assumed that $L_{p}\left(
S,\mu \right) $ is isometric to a separable space $L_{p}[0,1]$.

Let $E_{W}$ be the separable space of almost universal disposition from $%
\ast \ast \left( W^{f}\right) $. Since $\tau (W)=r$ and $r\leq p\leq 2$, $%
L_{p}\left( S,\mu \right) $ may be regarded as a subspace of $E_{W}$. Hence,
by the theorem 18, $w$ may be extended to an operator $\overline{w}:$ $%
L_{p}\left( S,\mu \right) \rightarrow G_{X}^{W}$, $\left\| \overline{w}%
\right\| \leq \lambda \left\| w\right\| $. Thus, $u$ has the factorization
\begin{equation*}
u=\overline{w}\circ i_{\mu ,p}\circ j,
\end{equation*}
and, hence, its $p$-integral norm $\iota _{p}\left( u\right) \leq \lambda
\pi _{p}\left( u\right) $.
\end{proof}

\begin{remark}
Analogously to [19] it may be shown that there exists a such Banach space $Y$
and a sequence of finite rank operators $u_{n}:Y\rightarrow G_{X}^{W}$ that
their $\pi _{p}$-norms are uniformly bounded and $\iota _{p}$-norms tend to
infinity with $n$:
\begin{equation*}
\sup_{n}\pi _{p}\left( u_{n}\right) \leq \lambda <\infty ;\text{ \ }%
\underset{n\rightarrow \infty }{\lim }\iota _{p}\left( u_{n}\right) =\infty
\end{equation*}
for any $p<r=\tau (W)$.

This implies that sub $B$-convex spaces $l_{1}\oplus _{1}W$ and $l_{1}\oplus
_{1}W_{1}$ from the same class $X^{f}$ and having different characteristics $%
\tau (W)\neq \tau (W_{1})$ generates different classes $\frak{Q}_{W}\left(
(l_{1}\oplus _{1}W)^{f}\right) $ and $\frak{Q}_{W}\left( (l_{1}\oplus
_{1}W_{1})^{f}\right) $ because corresponding spaces $G_{X}^{W}$ and $%
G_{X}^{W_{1}}$ are not isomorphic.
\end{remark}

\begin{remark}
It is obvious that the $l_{p}$-spectrum $S(\frak{Q}_{W}\left( (l_{1}\oplus
_{1}W)^{f}\right) =[1,q_{W}]$, where $q_{W}=\inf S(W)$.
\end{remark}

\section{J. Pisier like properties of spaces $G_{X}$}

Most results of this (and of the next) section are similar to those of J.
Pisier [1]. Their proofs are sufficiently short and are given for the
completeness.

\begin{theorem}
Let $X^{f}$ be $\frak{Q}$-closed class; $G_{X}$ be the corresponding
separable almost $\omega $-homogeneous space. Then its conjugate $\left(
G_{X}\right) ^{\ast }$is of cotype 2 and every operator $T:\left(
G_{X}\right) ^{\ast }\rightarrow l_{2}$ is 1-absolutely summing.
\end{theorem}

\begin{proof}
The result is a consequence of theorem 18 and the Pisier's result (cf. [1],
prop. 1.11).
\end{proof}

\begin{theorem}
Let $X^{f}$ be $\frak{Q}$-closed class; $G_{X}$ be the corresponding
separable almost $\omega $-homogeneous space. If $X^{f}$ is of cotype 2 then
any operator $v:G_{X}\rightarrow l_{2}$ is 1-absolutely summing.
\end{theorem}

\begin{proof}
Let $v:T:G_{X}\rightarrow l_{2}^{(n)}$. Then
\begin{equation*}
\pi _{1}\left( v\right) =\sup \{\pi \left( v\circ w\right) :\text{ \ }%
w:c_{0}\rightarrow G_{X};\text{ }\left\| w\right\| =1\}.
\end{equation*}

Since $G_{X}$ is of cotype 2, $\pi _{2}\left( w\right) \leq C\left\|
w\right\| =C$ and, according to the Pietsch theorem [20] $w$ can be factored
through $l_{2}$:
\begin{equation*}
w=w_{2}\circ w_{1};\text{ }w_{1}:c_{0}\rightarrow l_{2};\text{ \ }%
w_{2}:l_{2}\rightarrow G_{X};\text{ }\pi _{2}\left( w_{1}\right) \left\|
w_{2}\right\| \leq \pi _{2}\left( w\right) \leq C.
\end{equation*}

Since $\pi _{2}\left( \left( v\circ w_{2}\right) ^{\ast }\right) =\pi
_{2}\left( v\circ w_{2}\right) $,
\begin{equation*}
\pi _{2}\left( v\circ w_{2}\right) \leq \pi _{2}\left( w_{2}^{\ast }\right)
\left\| v\right\|
\end{equation*}
and
\begin{equation*}
\pi _{2}\left( w_{2}^{\ast }\right) \leq C\left\| w_{2}\right\| .
\end{equation*}

Since a composition of two 2-absolutely summing operators is 1-absolutely
summing one
\begin{eqnarray*}
\pi _{1}\left( v\circ w\right) &=&\pi _{1}\left( v\circ w_{2}\circ
w_{1}\right) \leq \pi _{2}\left( v\circ w_{2}\right) \pi _{2}\left(
w_{1}\right) \\
&\leq &C\left\| w_{2}\right\| \left\| v\right\| \pi _{2}\left( w_{1}\right)
\leq CC^{\prime }\left\| v\right\| .
\end{eqnarray*}

Thus, $\pi _{1}\left( v\right) \leq CC^{\prime }\left\| v\right\| $ and,
hence, $AS_{1}(G_{X},l_{2})=L\left( G_{X},l_{2}\right) $.
\end{proof}

\begin{remark}
The equality $AS_{1}(X,l_{2})=L\left( X,l_{2}\right) $ is equivalent to $%
AS_{2}(l_{2},X)=L\left( l_{2},X\right) $. Hence the following result is true.
\end{remark}

\begin{theorem}
Let $X^{f}$ be $\frak{Q}$-closed class; $G_{X}$ be the corresponding
separable almost $\omega $-universal space, which is of cotype 2. Then

\begin{itemize}
\item  $AS_{1}(G_{X},l_{2})=L\left( G_{X},l_{2}\right) $;

\item  $AS_{1}(\left( G_{X}\right) ^{\ast },l_{2})=L\left( \left(
G_{X}\right) ^{\ast },l_{2}\right) $;

\item  $AS_{2}(l_{2},G_{X})=L\left( l_{2},G_{X}\right) $;

\item  $AS_{2}(l_{2},\left( G_{X}\right) ^{\ast })=L\left( l_{2},\left(
G_{X}\right) ^{\ast }\right) $.
\end{itemize}
\end{theorem}

\begin{remark}
If the conditions of preceding theorem are satisfied, then

\begin{itemize}
\item  $AS_{1}(G_{X},Z)=AS_{2}\left( G_{X},Z\right) $;

\item  $AS_{1}(\left( G_{X}\right) ^{\ast },Z)=AS_{2}\left( \left(
G_{X}\right) ^{\ast },Z\right) $
\end{itemize}

for any arbitrary Banach space $Z$.

Indeed, if $u:G_{X}\rightarrow Z$ is 2-absolutely summing then it is
2-integral and, hence, can be factored through $l_{2}$. Since $%
AS_{1}(G_{X},l_{2})=L\left( G_{X},l_{2}\right) $, $u$ is 1-absolutely
summing. Same arguments prove the second equality.
\end{remark}

\begin{theorem}
Let $X^{f}$ be $\frak{Q}$-closed class; $G_{X}$ be the corresponding
separable almost $\omega $-homogeneous space, which is of cotype 2. There
exists a constant $\lambda <\infty $ such that for any finite rank operator $%
u:G_{X}\rightarrow G_{X}$ its nuclear norm
\begin{equation*}
\nu _{1}\left( u\right) \leq \left\| u\right\| .
\end{equation*}
\end{theorem}

\begin{proof}
According to [21], if $X^{\ast }$ and $Y$ both are of cotype 2 then there
exists such a constant $K<\infty $ that for any finite rank operator $%
v:Y\rightarrow X$%
\begin{equation*}
\nu _{1}\left( v\right) \leq K\inf \{\pi _{2}\left( w_{1}\right) ,\pi
_{2}\left( w_{2}\right) \},
\end{equation*}
where the infimum is taken over all possible factorizations $%
w_{1}:Y\rightarrow l_{2}$; $w_{2}:l_{2}\rightarrow X$; $w_{2}\circ w_{1}=v$.

Since $G_{X}$ and $\left( G_{X}\right) ^{\ast }$ both are of cotype 2, there
exists a constant $C<\infty $, which depends only on $G_{X}$, such that for
any finite rank operator $u:G_{X}\rightarrow G_{X}$ there exists its
factorization $u=u_{1}\circ u_{2}$ such that
\begin{equation*}
\inf \left\| u_{1}\right\| \left\| u_{2}\right\| \leq C\left\| v\right\| .
\end{equation*}

Since $\pi _{2}\left( u\right) \leq \pi _{1}\left( u\right) $ for any
operator $u$, there exists constants $C^{\prime }$ and $C^{\prime \prime }$%
such that
\begin{equation*}
\pi _{2}\left( u_{1}\right) \leq C^{\prime }\left\| u_{1}\right\| ;\text{ \ }%
\pi _{2}\left( \left( u_{2}\right) ^{\ast }\right) \leq C^{\prime \prime
}\left\| u_{2}\right\|
\end{equation*}
and, hence,
\begin{equation*}
\nu _{1}\left( v\right) \leq K\pi _{2}\left( u_{1}\right) \pi _{2}\left(
\left( u_{2}\right) ^{\ast }\right) \leq KC^{\prime }\left\| u_{1}\right\|
C^{\prime \prime }\left\| u_{2}\right\| \leq \lambda \left\| v\right\| .
\end{equation*}
\end{proof}

\begin{corollary}
Let $X^{f}$ be $\frak{Q}$-closed class; $G_{X}$ be the corresponding
separable almost $\omega $-homogeneous space, which is of cotype 2. Then

\begin{itemize}
\item  $AS_{1}\left( G_{X},l_{1}\right) =L\left( G_{X},l_{1}\right) $;

\item  $AS_{1}\left( \left( G_{X}\right) ^{\ast },l_{1}\right) =L\left(
\left( G_{X}\right) ^{\ast },l_{1}\right) $.
\end{itemize}
\end{corollary}

\begin{proof}
Since $G_{X}$ is an approximative envelope of a class $X^{f}$, it contains a
subspace, isomorphic to $l_{1}$. Let $u:G_{X}\rightarrow l_{1}$ be a finite
rank operator; $j:l_{1}\rightarrow G_{X}$ be an isometric embedding. Then
there exists such $\lambda <\infty $ that
\begin{equation*}
\pi _{1}\left( u\right) \leq \nu _{1}\left( j\circ u\right) \leq \lambda
\left\| j\circ u\right\| =\lambda \left\| u\right\| .
\end{equation*}

To prove the same result for $\left( G_{X}\right) ^{\ast }$ notice that $%
G_{X}$ contains a subspace isomorphic to $l_{2}$. So, $\left( G_{X}\right)
^{\ast \ast }$ contains an isomorphic image of $l_{2}$ too. Thus, there is a
surjection $h:\left( G_{X}\right) ^{\ast }\rightarrow l_{2}$. According to
theorem 22, $h$ is 1-absolutely summing. Since $h$ is not compact, from [22]
follows that $\left( G_{X}\right) ^{\ast }$ contains a subspace, isomorphic
to $l_{1}$.

Now we can repeat the proof of the first part of the corollary.
\end{proof}

\section{On the Grothendieck problem}

Let $X$, $Y$ be Banach spaces; $X\otimes Y$ be their \textit{algebraic
tensor product}. Define on $X\otimes Y$ two norms:
\begin{equation*}
\left\| u\right\| _{\wedge }=\inf \{\sum\nolimits_{i=1}^{n}\left\|
x_{i}\right\| \left\| y_{i}\right\| :u=\sum\nolimits_{i=1}^{n}x_{i}\otimes
y_{i}\},
\end{equation*}
where the infimum is taken over all possible representations of the tensor $%
u\in X\otimes Y$;
\begin{equation*}
\left\| u\right\| _{\vee }=\sup \{\left| \left\langle u,x^{\prime }\otimes
y^{\prime }\right\rangle \right| :x^{\prime }\in X^{\ast };\text{ }y^{\prime
}\in Y^{\ast };\text{ }\left\| x^{\prime }\right\| =\left\| y^{\prime
}\right\| =1\}.
\end{equation*}

The completition of $X\otimes Y$ with respect to the norm $\left\| \cdot
\right\| _{\wedge }$ is called \textit{the projective tensor product of }$X$%
\textit{\ and }$Y$ and is denoted by $X\overset{\wedge }{\otimes }Y$.

The completition of $X\otimes Y$ with respect to the norm $\left\| \cdot
\right\| _{\vee }$ is called \textit{the injective tensor product} \textit{%
of }$X$\textit{\ and }$Y$ and is denoted by $X\overset{\vee }{\otimes }Y$.

In [23], p.136 A. Grothendieck conceived that the coincidence of $X\overset{%
\wedge }{\otimes }Y$ and $X\overset{\vee }{\otimes }Y$ both in the algebraic
and in the topological sense would imply that at least one of these spaces
is of finite dimension.

In [1] J. Pisier disproved this conjecture by construction of a Banach space
$X_{P}$ with the property: $X_{P}\overset{\wedge }{\otimes }X_{P}\simeq X_{P}%
\overset{\vee }{\otimes }X_{P}$.

In this section it will be shown that all spaces of kind $G_{X}$, which are
of cotype 2, enjoy this property too.

The proof is based on ideas from [1], in particular, on a remark due to
S.V.~Kisliakov that is contained in [1] and will be used as the following
lemma.

\begin{lemma}
Let $L_{p}^{(n)}$ denotes the linear space $\mathbb{R}^{n}$, equipped with
the norm
\begin{equation*}
\left\| \left( x_{i}\right) _{i=1}^{n}\right\| _{p}=\left( \frac{1}{n}%
\sum\nolimits_{i=1}^{n}\left| x_{i}\right| ^{p}\right) ^{1/p}.
\end{equation*}

Then for every $n\in \mathbb{N}$ the space $L_{p}^{(n)}$ may be decomposed
into three parts $s_{n}^{1}$, $s_{n}^{2}$ and $s_{n}^{3}$: $%
L_{p}^{(n)}=s_{n}^{1}\oplus s_{n}^{2}\oplus s_{n}^{3}$ with the following
properties:

\begin{enumerate}
\item  The norms of $L_{2}^{(n)}$ and $L_{1}^{(n)}$ are uniformly equivalent
(with constants independent on $n$) on $s_{n}^{1}\oplus s_{n}^{2}$ and on $%
s_{n}^{2}\oplus s_{n}^{3}$.

\item  $dim(s_{n}^{i})$ tends to infinity with $n$; $\underset{n\rightarrow
\infty }{\lim }dim(s_{n}^{i})=\infty $ for $i=1,2,3$.
\end{enumerate}
\end{lemma}

Clearly, $L_{p}^{(n)}$ is isometric to $l_{p}^{(n)}$ and $%
s_{n}^{2}\hookrightarrow L_{p}^{(n)}/s_{n}^{3}$.

\begin{theorem}
Let $X^{f}$ be $\frak{Q}$-closed class; $G_{X}$ be the corresponding
separable almost $\omega $-homogeneous space, which is of cotype 2. Then
\begin{equation*}
G_{X}\overset{\wedge }{\otimes }G_{X}\simeq G_{X}\overset{\vee }{\otimes }%
G_{X}.
\end{equation*}
\end{theorem}

\begin{proof}
For convenience denote $G_{X}$ by $G$. Let $u\in G\otimes G$. Let us show
that norms $\left\| u\right\| _{\wedge }$ and $\left\| u\right\| _{\vee }$
are equivalent.

The tensor $u$ defines a finite rank operator $u:G^{\ast }\rightarrow G$.

Since $G$ is of cotype 2, there is a factorization $u=u_{2}\circ u_{1}$
through $l_{2}$; $\left\| u_{2}\right\| \left\| u_{1}\right\| \leq C\left\|
u\right\| $.

Since $u$ is weakly* continuous, it may be assumed that the same is true for
$u_{1}:G^{\ast }\rightarrow l_{2}$. Hence, $u_{1}=v^{\ast }$ for some finite
rank operator $v:l_{2}\rightarrow G$. Let $rg(v)=dim(vl_{2})=n$. Let $\left(
e_{i}\right) $ be a natural basis of $l_{2}$. It may be assumed that $%
v(e_{i})=0$ for $i>n$.

Let $l_{2}^{\left( m\right) }$ be realized as $L_{2}^{\left( m\right) }$.
Let $s_{m}^{2}\hookrightarrow L_{2}^{\left( m\right) }$ be as in lemma 1; $%
dim(s_{m}^{2})=n$. According the preceding results, it may be found an
operator $\overline{v}:L_{1}^{\left( m\right) }/s_{m}^{3}\rightarrow G$ such
that $\left\| \overline{v}\right\| \leq K\left\| v\right\| $ for some
constant $K<\infty $ that has a property:
\begin{equation*}
\overline{v}\mid _{s_{m}^{2}}=v\circ k^{-1}\mid _{s_{m}^{3}};\text{ }%
k:l_{2}^{\left( n\right) }\rightarrow s_{m}^{2};\text{ }\left\| k\right\|
\left\| k^{-1}\right\| =1.
\end{equation*}

Clearly, since $s_{m}^{2}\hookrightarrow L_{1}^{\left( m\right) }/s_{m}^{3}$
, $k$ may be regarded as an isometric embedding of $l_{2}^{\left( n\right) }$
into $L_{1}^{\left( m\right) }/s_{m}^{3}$.

From the property 1 of lemma 1 of the decomposition $L_{p}^{(m)}=s_{m}^{1}%
\oplus s_{m}^{2}\oplus s_{m}^{3}$ follows that there exists a finite rank
operator $T_{n}:l_{1}/s_{m}^{3}\rightarrow l_{1}/s_{m}^{3}$ which is
identical on $ks_{m}^{2}$, which vanished on its orthogonal complement in $%
l_{2}$and whose norm is bounded by the constant $C^{\prime }$ that does not
depend on $n$: $\left\| T_{n}\right\| \leq C^{\prime }$.

Clearly, $v=\overline{v}\circ T_{n}\circ k$. Hence, $u_{1}=k^{\ast }\circ w$%
, where $w=\left( \overline{v}\circ T_{n}\right) ^{\ast }$;
\begin{equation*}
\left\| w\right\| \leq C^{\prime }\left\| \overline{v}\right\| \leq
C^{\prime }K\left\| v\right\| =KC^{\prime }\left\| u_{1}\right\| .
\end{equation*}

Notice that $k^{\ast }:L_{\infty }^{\left( m\right) }\rightarrow
l_{2}^{\left( n\right) }$ has a factorization $k^{\ast }=P\circ J$, where
the map $J:L_{\infty }^{\left( m\right) }\rightarrow L_{1}^{\left( m\right)
} $ is the natural injection and $P:L_{1}^{\left( m\right) }\rightarrow
l_{2}^{\left( n\right) }$.

Therefore $u_{1}=P\circ J\circ w$ and, thus, $u=u_{2}\circ P\circ J\circ w$.

Consider $T=u_{2}\circ P:l_{1}^{\left( m\right) }\rightarrow G$. From the
theorem 18 follows that $T$ has an extension $\overline{T}:L_{1}\rightarrow
G $ with $\left\| \overline{T}\right\| \leq K^{\prime }\left\| T\right\| $
for some $K^{\prime }<\infty $; $\overline{T}\mid _{L_{1}^{\left( m\right)
}}=T$.

Let $\overline{J}:L_{\infty }\rightarrow L_{1}$ and $i:L_{\infty }^{\left(
m\right) }\rightarrow L_{\infty }$ be natural injections. Then
\begin{equation*}
u=\overline{T}\circ \overline{J}\circ i\circ w.
\end{equation*}

The operator $\overline{T}\circ \overline{J}\circ i$ is 1-integral. Its
1-integral norm is estimated as
\begin{equation*}
\iota _{1}\left( \overline{T}\circ \overline{J}\circ i\right) \leq \left\|
\overline{T}\right\| \left\| i\right\| =\left\| \overline{T}\right\| \leq
K^{\prime }\left\| u_{2}\right\| .
\end{equation*}

Since $w$ is a finite rank operator, the space $E=wG$ is of finite
dimension. $w$ may be represented as $h\circ \overline{w}$, where $%
h:E\rightarrow L_{\infty }^{\left( m\right) }$ is the identical embedding
and $\overline{w}:G\rightarrow E$ is the restriction of $w$. Since $w$ is
weakly* continuous, the same is true for $\overline{w}$.

Consider the operator $U=\overline{T}\circ \overline{J}\circ i\circ
h:E\rightarrow G$.

It was shown that $\iota _{1}\left( U\right) \leq K^{\prime }\left\|
u_{2}\right\| $. Since $E$ is of finite dimension, the 1-integral norm $%
\iota _{1}\left( U\right) $ of $U$ is equal to its 1-nuclear norm $\nu
_{1}\left( U\right) $.

Hence, $U$ may be identified with a tensor in $E^{\ast }\otimes G$; $\left\|
U\right\| _{\wedge }=\iota _{1}\left( U\right) \leq K^{\prime }\left\|
u_{2}\right\| $.

Since $\overline{w}$ is weakly* continuous, the tensor $u=U\circ v\in
G\otimes G$ satisfies:
\begin{equation*}
\left\| u\right\| _{\wedge }\leq \left\| U\right\| _{\wedge }\left\|
w\right\| _{\vee }\leq K^{\prime }\left\| u_{2}\right\| \left\| w\right\|
_{\vee }\leq KK^{\prime }C^{\prime }\left\| u_{2}\right\| \left\|
u_{1}\right\| \leq KK^{\prime }CC^{\prime }\left\| u\right\| .
\end{equation*}

Thus, norms $\left\| \cdot \right\| _{\wedge }$ and $\left\| \cdot \right\|
_{\vee }$ are equivalent on $G\otimes G$.
\end{proof}

\section{Projection constants}

The following result is similar to [1], remark 3.4.

\begin{theorem}
Let $X^{f}$ be $\frak{Q}$-closed class; $G_{X}$ be the corresponding
separable almost $\omega $-homogeneous space, which is of cotype 2. Let $%
T:G_{X}\rightarrow G_{X}$ be a finite rank operator, which fixes a finite
dimensional subspace $A\hookrightarrow G_{X}$. Then
\begin{equation*}
\left\| T\right\| \geq C\sqrt{dim(A)},
\end{equation*}
where $C$ is a constant, which depends only on $G_{X}$.
\end{theorem}

\begin{proof}
Since $G_{X}$ and $\left( G_{X}\right) ^{\ast }$ both are of cotype 2, there
exists a constant $C<\infty $ such that for every finite rank operator $%
T:G_{X}\rightarrow G_{X}$
\begin{equation*}
\inf \{\left\| v\right\| \left\| w\right\| :T=v\circ w;\text{ }%
v:l_{2}\rightarrow G_{X};\text{ }w:G_{X}\rightarrow l_{2}\}\leq C\left\|
T\right\| .
\end{equation*}

Thus, since $AS_{1}(G_{X},l_{2})=L\left( G_{X},l_{2}\right) $,
\begin{equation*}
\pi _{2}\left( T\right) \leq C^{\prime }\left\| v\right\| \left\| w\right\|
\leq CC^{\prime }\left\| T\right\| .
\end{equation*}

Let $\left( \lambda _{i}\right) _{i\geq 1}$ be eigenvalues of $T$. According
to [24], 27.4.6,
\begin{equation*}
\sqrt{\sum \left| \lambda _{i}\right| ^{2}}\leq \pi _{2}\left( T\right) \leq
CC^{\prime }\left\| T\right\| .
\end{equation*}

If $T$ fixes $A\hookrightarrow G_{X}$ and $dim(A)=n$, then there exists at
least $n$ of eigenvalues $\left( \lambda _{i}\right) _{i\geq 1}$ of $\left|
\lambda _{i}\right| =1$. Thus, $\left\| T\right\| \geq C\sqrt{n}$.
\end{proof}

\begin{corollary}
If $G_{X}$ is of cotype 2 then there exists a constant $C<\infty $ such that
for any projection $P:G_{X}\rightarrow G_{X}$
\begin{equation*}
\left\| P\right\| \geq C\left( rg(P)\right) ^{1/2}.
\end{equation*}
\end{corollary}

It may be presented a different proof of this corollary which extends it on
spaces of kind $G_{X}$ of arbitrary cotype $q<\infty $. Recall some
definitions.

Let $X$ be a Banach space; $Y\hookrightarrow X.$A \textit{projection constant%
} $\lambda (A\hookrightarrow X)$ is given by
\begin{equation*}
\lambda (A\hookrightarrow X)=\inf \{\left\| P\right\| :\text{ }%
P:X\rightarrow A\},P^{2}=P,
\end{equation*}
i.e. $P$ runs all projections of $E_{X}$ onto $A$.

\textit{Relative projection constants} $\lambda (A,X)$ and $\lambda (A,%
\mathcal{K})$, where $\mathcal{K}$ is a class of Banach spaces are defined
as follows:
\begin{equation*}
\lambda (A,X)=\sup \{\lambda (iA\hookrightarrow X):\text{ \ }i:A\rightarrow
X\},
\end{equation*}
where $i$ runs all isometric embeddings of $A$ into $X$;
\begin{equation*}
\lambda (A,\mathcal{K})=\sup \{\lambda (A,X):X\in \mathcal{K}\}.
\end{equation*}

The \textit{absolute projection constant }$\lambda (A)=\lambda (A,\mathcal{B}%
)$, where $\mathcal{B}$ is the class of all Banach spaces.

\begin{theorem}
Let $G$ be a separable almost $\omega $-homogeneous space; $A\hookrightarrow
G$. Then
\begin{equation*}
\lambda (A\hookrightarrow G)=\lambda (iA\hookrightarrow G)=\lambda (A,G)
\end{equation*}
where $i$ is an arbitrary isometric embedding.
\end{theorem}

\begin{proof}
Obviously, any almost $\omega $-homogeneous space has the property:

\begin{itemize}
\item  \textit{For any pair }$A$\textit{, }$B$\textit{\ of isometric
subspaces of }$G$\textit{\ under an isometry }$j:A\rightarrow B$\textit{\
and every }$\varepsilon >0$\textit{\ there exists an isomorphical
automorphism }$u:G\rightarrow G$\textit{\ such that }$u\mid _{A}=j$\textit{\
and }$\left\| u\right\| \left\| u^{-1}\right\| \leq 1+\varepsilon $\textit{.}
\end{itemize}

Clearly, this property implies the theorem.
\end{proof}

\begin{theorem}
Let $X^{f}$ be $\frak{Q}$-closed class; $G_{X}$ be the corresponding
separable almost $\omega $-homogeneous space. Let $G_{X}$ be of cotype $q$.
Then there exists a constant $C$ such that for any projection $%
P:G_{X}\rightarrow G_{X}$ of finite range $rg(P)=n$
\begin{equation*}
\left\| P\right\| \geq Cn^{1/q}.
\end{equation*}
\end{theorem}

\begin{proof}
Assume at first that $PG_{X}=l_{2}^{(n)}$. Since $L_{1}[0,1]$ is almost $%
\omega $-homogeneous (this follows from the amalgamation property of $\frak{M%
}(\left( l_{1}\right) ^{f})$), $\lambda (il_{2}^{(n)}\hookrightarrow
L_{1}[0,1])$ does not depend on the isometric embedding $i$. At the same
time it may be assumed that $L_{1}[0,1]\hookrightarrow G_{X}$. Since $G_{X}$
is almost $\omega $-homogeneous, $\lambda (jl_{2}^{(n)}\hookrightarrow
G_{X}) $ also does not depend on $j$. Thus,
\begin{equation*}
\left\| P\right\| \geq \lambda (l_{2}^{(n)},G_{X})\geq \lambda
(l_{2}^{(n)},L_{1}[0,1]).
\end{equation*}

However, it is known that $\lambda (l_{2}^{(n)},L_{1}[0,1])\asymp \sqrt{n}$
for large $n$ (up to constant).

The same estimate is true when $PG_{X}$ is $\left( 1+\varepsilon \right) $%
-close to $l_{2}^{(n)}$, i.e. when
\begin{equation*}
d\left( l_{2}^{(n)},PG_{X}\right) \leq \left( 1+\varepsilon \right) .
\end{equation*}

Let $PG_{X}=A$; $dim(A)=2n$. Clearly, $A$ may be represented as $%
A=l_{1}^{(k)}\oplus _{1}A^{\prime }$, where $dim(A^{\prime })=2n-k$ and for
some $p>1$ the $p$-type constant $T_{p}(A^{\prime })\leq \lambda $, where $%
\lambda $ depends only on $G_{X}$.

Consider two cases.

1. $n\leq k$. Since $G_{X}$ is of cotype $q$ then its $q$-cotype constant $%
C_{q}\left( G_{X}\right) =C^{\prime }<\infty $. Hence, the same is true for $%
A^{\prime }$: $C_{q}\left( A^{\prime }\right) \leq C^{\prime }$.

According to the Dvoretzki theorem, $A^{\prime }$ contains a subspace that
is sufficiently close to $l_{2}^{(m)}$ for some $m<n$. More precise
estimates were given in [25]. Namely, in particular there was shown that a
dimension $m$ of a subspace $B$ of $A^{\prime }$, that is $\left(
1+\varepsilon \right) $-close to $l_{2}^{(m)}$ may be estimated as
\begin{equation*}
m=dim(B)\geq c(C^{\prime },\varepsilon ,\lambda )n^{2/q},
\end{equation*}
where $c(C^{\prime },\varepsilon )$ is a constant which depends only on $%
C^{\prime }$, $\varepsilon $ and $\lambda $, where the type $p$-constant of $%
A^{\prime }$ is bounded by $\lambda $ for some $p>1$:
\begin{equation*}
\exists \left( p>1\right) :T_{p}(A^{\prime })\leq \lambda .
\end{equation*}

For any operator $F:G_{X}\rightarrow G_{X}$ of finite range, which fixes $B$%
, $\left\| F\right\| \geq \lambda (B,G_{X})$. Hence, if $P:G_{X}\rightarrow
A $ and $P^{2}=P$,
\begin{equation*}
\left\| P\right\| \geq C^{\prime \prime }\left( c(C^{\prime },\varepsilon
)n^{2/q}\right) ^{1/2}\geq Cn^{1/q}.
\end{equation*}

2. $k\leq n$. Let $P:G\rightarrow l_{1}^{(n)}$ be a projection.

The conjugate operator $P^{\ast }:c_{0}^{\left( n\right) }\rightarrow \left(
G_{X}\right) ^{\ast }$ is an operator of isomorphic embedding.

Since $\left( G_{X}\right) ^{\ast }$ is of cotype 2, $P^{\ast }$ is
2-absolutely summing. Hence, $\left\| P^{\ast }\right\| \geq c\sqrt{n}$ for
a some constant $c$ that depends only on $G_{X}$. Since $\left\| P^{\ast
}\right\| =\left\| P\right\| $, the theorem is proved.
\end{proof}

\begin{corollary}
If $W$ is a $B$-convex Banach space of cotype $q$ then there exists a
constant $C<\infty $ such that for any its finite dimensional subspace $A$%
\begin{equation*}
\lambda (A)\geq Cn^{1/q}
\end{equation*}
\end{corollary}

\section{Factorization of operators}

Grothendieck asked ([23], problem 2, p. 78): \textit{whether there exists an
absolutely summing operator which is not }$L_{1}$\textit{-factorable}?

Gordon and Lewis [26] answered this question in affirmative. At the same
time they showed that there exists a wide class of Banach spaces $X$with the
property:

(GL)\ \textit{Every absolutely summing operator }$u:X\rightarrow Y$\textit{\
where }$Y$\textit{\ is an arbitrary Banach space is }$L_{1}$\textit{%
-factorable. }

This class consists of \textit{spaces with the local unconditional structure}
(shortly, l.u.st.) which were introduced in [26]. Recall that, according to
[27], a Banach space $X$ has the l.u.st. if and only if its second conjugate
is isomorphic to a complemented subspace of some Banach lattice.
Particularly, the class of all Banach spaces having the l.u.st. contains all
Banach lattices. The result of [26] further was improved (see e.g. [28],
p.107); now it is known that if a Banach space $X$ is of cotype 2 and has
the l.u.st. then each its subspace enjoy the property (GL).

Gordon and Lewis [26] defined for a given Banach space $X$ a constant
\begin{equation*}
GL\left( X\right) =\sup \{\gamma _{1}\left( u\right) /\pi _{1}\left(
u\right) :u:X\rightarrow Y\},
\end{equation*}
where the supremum is taken over all Banach spaces $Y$ and all absolutely
summing operators $u:X\rightarrow Y$ (here and below $\pi _{1}\left(
u\right) $denotes the absolutely summing norm of an operator $u$). In their
solution of the aforementioned Grothendieck problem a sequence $\left(
u_{n}\right) _{n<\infty }$of finite rank operators ($u_{n}:X\rightarrow
l_{2}^{(m_{n})}$) was constructed in such a way that $\sup \{\gamma
_{1}\left( u_{n}\right) /\pi _{1}\left( u_{n}\right) :n<\infty \}=\infty .$

Let us define a\textit{\ finite dimensional version of the Gordon-Lewis
constant:}
\begin{equation*}
GL_{fin}(X)=\sup\{\gamma_{1}\left( u\right) /\pi_{1}\left( u\right) :\text{
\ }u:X\rightarrow l_{2};\text{ \ }\dim u(X)<\infty\}.
\end{equation*}

It is clear that $GL_{fin}(X)=\infty $ implies that $GL(X)=\infty $.

However in a general case may be a situation when for a given Banach space $%
X $ constants $GL_{fin}(X)$ and $GL(X)$ are essentially different. To be
more exact, it will be shown that for Banach spaces of kind $G_{X}$ the
constant $GL_{fin}(G_{X})$ is finite and $GL(G_{X})$ is infinite.

\begin{theorem}
Let $X^{f}$ be a $\frak{Q}$-closed class of cotype 2; $G_{X}$ be the
corresponding separable almost $\omega $-homogeneous space. Then $GL\left(
G_{X}\right) =\infty $ and $GL_{fin}(G_{X})<\infty $.
\end{theorem}

\begin{proof}
Certainly, $G_{X}$ has subspaces that are isometric to $l_{1}$ and to $l_{2}$%
. Let $l_{2}$ be isometric to a subspace $X_{2}\hookrightarrow G_{X}$ and $%
i:l_{2}\rightarrow X_{2}$ be the corresponding isometry. Since $G_{X}$ is of
cotype 2, each operator $u\in B(G_{X},l_{2})$ is absolutely summing . Let $%
u:G_{X}\rightarrow l_{2}$. The composition $w=i\circ u$ acts on $G_{X}$, is
absolutely summing and is $L_{1}$-factorable if and only if $u$ is. Assume
that $u$ has a finite dimensional range: $dim\left( u(G_{X})\right) <\infty $%
. Then $w$ is also a finite rank operator. Hence, its nuclear norm $\nu
_{1}\left( w\right) $ may be estimated as $\nu _{1}\left( w\right) \leq
\lambda \left\| w\right\| $, where $\lambda $ does not depend on $w$. Thus,
\begin{equation*}
\gamma _{1}\left( u\right) =\gamma _{1}\left( w\right) \leq \iota _{1}\left(
w\right) \leq \nu _{1}\left( w\right) \leq \lambda \left\| w\right\| \leq
\lambda \pi _{1}\left( w\right) .
\end{equation*}

This sequence of inequalities implies that
\begin{equation*}
\gamma _{1}\left( u\right) \leq \lambda \pi _{1}\left( u\right) .
\end{equation*}

Recall that $u$ was an arbitrary operator of finite dimensional range.
Hence, $GL_{fin}(G_{X})<\infty $.

Let us show that for a special choosing of $u:G_{X}\rightarrow l_{2}$ such
inequalities cannot be longer true. Indeed, since $G_{X}$ contains a
subspace isomorphic to $l_{1}$ then its conjugate $X^{\ast }$ contains a
subspace (say, $Y_{2}$), which is isomorphic to $l_{2}$. Let $%
j:l_{2}\rightarrow Y_{2}$ be the corresponding isomorphism; $\left\|
j\right\| \left\| j^{-1}\right\| <\infty $. Because of reflexivity, $Y_{2}$
is weakly* closed and, hence, there exists a quotient map $h:X\rightarrow
(Y_{2})^{\ast }$ which generates a surjection $H:X\rightarrow l_{2}$: $%
H=(j_{2})^{\ast }\circ h$. From properties of $G_{X}$ follows that the
conjugate operator $H^{\ast }$ is $L_{\infty }$-factorable and is
2-absolutely summing. Clearly, $H^{\ast }$ induces an isomorphism between $%
l_{2}$ and $H^{\ast }(l_{2})$. However this property contradicts with $\pi
_{2}\left( H^{\ast }\right) <\infty $. Hence, $H$ cannot be factored through
$L_{1}\left( \mu \right) $.
\end{proof}

\begin{definition}
Let $Y$ be a Banach space of cotype 2. It will be said that $X\in \mathcal{B}
$ is the \textit{Pisier's space for }$Y$, if following conditions are
satisfy.

\begin{enumerate}
\item  $Y$ is isometric to a subspace of $X.X$ and $X^{\ast }$ both are of
cotype 2.

\item  Any operator $u:X\rightarrow l_{2}$ is absolutely summing; any
operator $v:X^{\ast }\rightarrow l_{2}$ is absolutely summing.

\item  There is a constant $\lambda _{X}$ such that for any finite rank
operator $u:X\rightarrow X$ its nuclear norm $\nu _{1}(u)$ satisfies the
estimate $\nu _{1}(u)\leq \lambda _{X}\left\| u\right\| $.

\item  If $X$ contains a subspace isomorphic to $l_{1}$ then its conjugate $%
X^{\ast }$ contains a subspace isomorphic to $l_{2}$.

\item  Any operator $w:l_{2}\rightarrow X$ is 2-absolutely summing; the same
is true for every operator $s:l_{2}\rightarrow X^{\ast }$.
\end{enumerate}
\end{definition}

The famous result of J. Pisier [1] shows the existence of the corresponding
Pisier's space $X$ for any Banach space $Y$ of cotype 2.

From preceding sections follows that for any $B$-convex Banach space $W$ of
cotype 2 the space $G_{X}$ - the almost $\omega $-homogeneous space for a
class $\frak{Q}_{W}\left( (l_{1}\oplus _{1}W)^{f}\right) $ is the Pisier's
space for\textit{\ }$W$.

\begin{problem}
\textit{Whether any }$c$\textit{-convex Banach space }$X$\textit{\ generates
a sub }$B$\textit{-convex class }$X^{f}$\textit{? }
\end{problem}

If this question has an affirmative answer then for any Banach space $Z$ of
cotype 2 the corresponding Pisier's space contains among spaces of kind $%
G_{X}$.

\begin{theorem}
For every Banach space of cotype 2, $GL_{fin}(X)$ is finite.
\end{theorem}

\begin{proof}
Let a Banach space $Y$ is of cotype 2. It can be isometricaly embedded in a
corresponding Pisier's space $X_{Y}$, which has the same properties as was
listed above. Let $T:Y\rightarrow l_{2}$ be an absolutely summing operator.
Then $T$ is 2-absolutely summing and, consequently, 2-integral, i.e. admits
a factorizations $T=k\circ j\circ s$, where $s:X\rightarrow L_{\infty }(\mu )
$; $k\circ j:L_{\infty }(\mu )\rightarrow l_{2}$. Because of injectivity of $%
L_{\infty }(\mu )$, $s$ can be extended to $\overset{\sim }{s}%
:X_{Y}\rightarrow L_{\infty }(\mu )$. The corresponding operator $\overset{%
\sim }{T}=k\circ j\circ \overset{\sim }{s}:X_{Y}\rightarrow l_{2}$ is
2-absolutely summing too and extends $T$. Recall that from properties of $%
X_{Y}$ follows that $\overset{\sim }{T}$ is absolutely summing. Assume that $%
T$ has a finite dimensional range. The same is true for $\overset{\sim }{T}$%
. From theorem 1 (that uses only listed above properties of Pisier's spaces)
follows that $\overset{\sim }{T}$ is $L_{1}$-factorable. Hence $T$ is also $%
L_{1}$-factorable.
\end{proof}

\begin{remark}
From [29] it follows that if there exists some $p>2$ such that the space $%
l_{p}$ is finitely representable in $X\in \mathcal{B}$ then $X$ contains a
subspace $Z$ for which $GL_{fin}(Z)$ is infinite.
\end{remark}

\section{References}

\begin{enumerate}
\item  Pisier G. \textit{Counterexample to a conjecture of Grothendieck},
Acta Math. \textbf{151} (1983) 181-208

\item  Schwartz L. \textit{Geometry and probability in Banach spaces}, Bull.
AMS \textbf{4:2} (1981) 135-141

\item  Maurey B., Pisier G. \textit{S\'{e}ries de variables al\'{e}atoires
vectori\'{e}lles ind\'{e}pendantes et propri\'{e}t\'{e}s g\'{e}om\'{e}%
triques des espaces de Banach}, Studia Math. \textbf{58} (1976) 45-90

\item  J\'{o}nsson B. \textit{Universal relational systems}, Math. Scand.
\textbf{4} (1956) 193-208

\item  Positselski E.D., Tokarev E.V. \textit{Amalgamation of classes of
Banach spaces }(in Russian), \textbf{XI} All-Union School on Operator Theory
in Functional Spaces. A book of abstracts, Tcherlyabinsk (1986) 89

\item  Zalgaller V.A., Reshetniak Yu.G. \textit{On straightened curves,
additive vector-functions and merging of segments }(in Russian), Vestnik LGU
\textbf{2} (1954) 45-67

\item  Rudin W. $L_{p}$\textit{-isometries and equimeasurability,} Indiana
Math. J. \textbf{25:3} (1976) 215-228

\item  Lusky W. \textit{Some consequences of W. Rudin's paper ''}$L_{p}$%
\textit{\ - isometries and equimeasurability''}, Indiana Univ. Math. J.
\textbf{27:5} (1978) 859-866

\item  Linde W. \textit{Moments of measures in Banach spaces}, Math. Ann.
\textbf{258} (1982) 277-287

\item  Gurarii V.I.\textit{\ Spaces of universal disposition, isotropic
spaces and the Mazur problem on rotations in Banach spaces}, Sibirsk. Mat.
Journ. (in Russian) \textbf{7} (1966) 1002-1013

\item  Pisier G. \textit{On the duality between type and cotype}, Lect.
Notes in Math. \textbf{939} (1982)

\item  Stern J. \textit{Ultrapowers and local properties in Banach spaces},
Trans. AMS \textbf{240} (1978) 231-252

\item  Tokarev E.V. \textit{Injective Banach spaces in the finite
equivalence classes }(transl. from Russian), Ukrainian Mathematical Journal
\textbf{39:6} (1987) 614-619

\item  Gelfand I.M. \textit{Abstracte Functionen und lineare Operatoren},
Mat. Sbornik, \textbf{4} (1938) 235-284

\item  Odell E., Rosenthal H.P. \textit{A double - dual characterization of
separable Banach spaces, containing }$\mathit{l}_{\mathit{1}}$, Israel. J.
Math. \textbf{20:3-4} (1975) 375-384

\item  Szlenk W. \textit{The non - existence of a separable reflexive Banach
space universal for all separable reflexive Banach spaces}, Studia Math.
\textbf{30} (1968) 53-61

\item  Wojtaszczyk P. \textit{On separable Banach spaces, containing all
separable reflexive Banach spaces,} Studia Math\textit{. }\textbf{37} (1970)
197-202

\item  Pe\l czy\'{n}ski A. \textit{On the impossibility of embedding of the
space }$L$\textit{\ in certain Banach spaces}, Colloq. Math. \textbf{8}
(1961) 199-203

\item  Pe\l czy\'{n}ski A\textit{. }$p$\textit{-integral operators commuting
with group representations and examples of quasi }$p$\textit{-integral
operators, which are not }$p$\textit{-integral}, Studia Math. \textbf{33}
(1969) 63-70

\item  Pietsch A. \textit{Absolut }$p$\textit{-summierende Abbildungen in
normierten Raumen}, Studia Math. \textbf{28} (1967) 333-353

\item  Maurey B. \textit{Un th\'{e}or\`{e}me de prolongement}, C. R. Acad.
Sci. Paris, S\'{e}r.A, \textbf{24} (1974) 329-332

\item  Rosenthal H.P. \textit{Pointwise compact subsets of the first Baire
class}, Amer. J. Math. \textbf{99} (1977) 362-378

\item  Grothendieck A. \textit{R\'{e}sum\'{e} de la th\'{e}orie m\'{e}trique
des produits tensoriels} \textit{topolo-giques,} Bol. Soc. Math. Sa\~{o}
Paulo, \textbf{8} (1956), 1-79

\item  Pietsch A. \textit{Operator ideals}, Amsterdam: North-Holland, 1980

\item  Figiel T., Lindenstrauss J., Milman V. \textit{The dimension of
almost spherical sections of convex bodies}, Acta Math. \textbf{139} (1977)
53-94

\item  Gordon Y., Lewis D.R. \textit{Absolutely summing operators and local}
\textit{unconditional structures}, Acta Math., \textbf{133 }(1974), 27-48

\item  Figiel T., Johnson W.B.,Tzafriri L. \textit{On Banach lattices and
spaces having local unconditional structure with applications to Lorentz
function spaces}, J. Approx. Theory, \textbf{13:4} (1975), 392-412

\item  Pisier J. \textit{Factorization of linear operators and geometry of
Banach spaces}, Publ.: Conf. Board of Math. Sci. by Amer. Math. Soc.,
Providence, Rhode Island, 1986, 153 p.

\item  Figiel T., Kwapie\'{n} S., Pe\l czy\'{n}ski A. \textit{Sharp
estimates for the constants of local unconditional structure of Minkowski
spaces}, Bull. Acad. Pol. Sci. S\'{e}r. Math. Astron. et Phys. \textbf{25 }%
(1977) 1221 - 1226
\end{enumerate}

\end{document}